\documentstyle[amssymb,amsfonts,12pt]{amsart}
\newtheorem{theorem}{Theorem}[section]
\newtheorem{lemma}[theorem]{Lemma}
\newtheorem{corollary}[theorem]{Corollary}
\newtheorem{proposition}[theorem]{Proposition}
\newtheorem{definition}[theorem]{Definition}

\theoremstyle{remark}
\newtheorem{remark}[theorem]{Remark}

\def\Z{{\bf Z}}
\def\R{{\bf R}}
\def\U{{\bf U}}
\def\C{{\bf C}}
\def\N{{\bf N}}
\def\E{{\mbox{\rm I\kern-.22em E}}}
\def\P{{\bf P}}
\def\D{{\bf D}}
\def\T{{\bf T}}
\def\S{{\bf S}}

\def\size{{\rm size}}

\def\K{{\mathcal K}}
\def\\P{{\mathcal P}}

\def\F{{\mathcal F}}

\def\G{{\mathcal G}}
\def\B{{\mathcal S}}
\def\I{{\mathcal I}}

\def\dist{{\rm dist}}

\def\BMO{{\operatorname{BMO}}}

\def\bas{\begin{align*}}
\def\eas{\end{align*}}
\def\bi{\begin{itemize}}
\def\ei{\end{itemize}}
\newenvironment{proof}{\noindent {\bf Proof} }{\endprf\par}
\def \endprf{\hfill  {\vrule height6pt width6pt depth0pt}\medskip}
\def\emph#1{{\it #1}}

\begin{document}
\title{Pointwise convergence of the ergodic bilinear Hilbert transform}

\author{Ciprian Demeter}
\address{Department of Mathematics, UCLA, Los Angeles CA 90095-1555}
\email{demeter@@math.ucla.edu}

\begin{abstract}
Let ${\bf X}=(X, \Sigma, m, \tau)$ be a dynamical system. We prove that the bilinear series $\sideset{}{'}\sum_{n=-N}^{N}\frac{f(\tau^nx)g(\tau^{-n}x)}{n}$ converges almost everywhere for each $f,g\in L^{\infty}(X).$ We also give a proof along the same lines of Bourgain's analog result for averages.
\end{abstract}
\maketitle

\section{Introduction}
\label{sec:intro}

Let ${\bf X}=(X, \Sigma, m, \tau)$ be a dynamical system, i.e. a complete probability space $(X,\Sigma,m)$ endowed with an invertible bimeasurable transformation $\tau:X\to X$ such that $m{\tau^{-1}}=m$. The starting point in this discussion is a result proved by Bourgain for bilinear averages.

\begin{theorem}[\cite{Bo10}]
\label{pontconvave}
For each $f,g\in L^{\infty}(X)$ the averages
\begin{equation}
\label{e.ave.Bou}
\frac1N\sum_{n=0}^{N-1}f(\tau^nx)g(\tau^{-n}x)
\end{equation}
converge for almost every $x$.
\end{theorem}

Bourgain's method consists of turning the issue of almost everywhere  convergence into a quantitative problem regarding multipliers on the torus, which are investigated by using classical Fourier analysis. An important reduction of ergodic theoretic nature in his argument concerns the fact that $g$ can be assumed to be orthogonal to the linear space $L^2({\bf K})$, where ${\bf K}=(X, \K, m)$ and $\K\subset\Sigma$ is the $\sigma$- algebra generated by the eigenfunctions of $\tau$. This is because the convergence is trivial in the case $g$ is an eigenfunction, as it is easily seen from Birkhoff's ergodic theorem \cite{Bi}. This reduction has the following  consequence for the spectral behavior of $g$
\begin{equation}
\label{WWofBou}
\lim_{N\to\infty}\sup_{|z|=1}|\frac1N\sum_{n=0}^{N-1}g(\tau^{-n}x)|=0,
\end{equation}
see \cite{Asa} for a proof of this and of some related results. Using this, Bourgain  identifies the limit to be 0 for such a $g$. More  generally  
$$\lim_{N\to\infty}\frac1N\sum_{n=0}^{N-1}f(\tau^nx)g(\tau^{-n}x)=\lim_{N\to\infty}\frac1N\sum_{n=0}^{N-1}f(\tau^nx)P_{{\bf K}}g(\tau^{-n}x),$$
for each $g\in L^{\infty}(X)$, where $P_{{\bf K}}g$ is the projection of $g$ onto $L^2({\bf K})$. A different way of putting  this is to say that the Kronecker factor ${\bf K}$ is a characteristic factor for the almost everywhere convergence of the averages ~\eqref{e.ave.Bou}.

In this paper we will prove the convergence of the ergodic bilinear Hilbert transform.

\begin{theorem}
\label{pontconvser}
For each $f,g\in L^{\infty}(X)$ the series
$$\sideset{}{'}\sum_{n=-N}^{N}\frac{f(\tau^nx)g(\tau^{-n}x)}{n}$$
converges for almost every $x$.
\end{theorem}
\begin{remark}
As a consequence of the above and of the bilinear maximal inequality for very general kernels in \cite{La2}, it follows that  Theorem ~\ref{pontconvser}
holds for all $f\in L^{p}(X)$, $g\in L^{q}(X)$, whenever $1<p,q\le \infty$ and $\frac{1}p+\frac1{q}<\frac32$.
\end{remark}

Bourgain's approach does not seem to be applicable to the context of series, in part due to the fact that the characteristic factors for weighted operators other than the usual averages are much less understood, and probably of less relevance to the essence of the problem. In particular, ~\eqref{WWofBou} fails for the series above, and the Kronecker factor seems to be of no immediate relevance to the problem. 

We prove Theorem ~\ref{pontconvser} using time-frequency harmonic analysis, and by a similar argument we also give a new proof of Theorem ~\ref{pontconvave}. Our methods  will not perceive the difference between the differentiation and the singular integral versions of the above, due to a common decomposition of both operators into discrete model sums.

Interestingly, our argument does not appeal to characteristic factors or in general to any concrete spectral analysis. Moreover, only little  ergodic theory is needed in the whole argument, when integration along individual orbits allows us to transfer certain oscillation inequalities from harmonic analysis. However, the structure of the Kronecker factor is deeply rooted into our approach. Since the (linear) exponentials $e^{i\lambda x}$ are the eigenfunctions for rotations on the torus, it is probably the case  that their  presence  in the wave packet decomposition of $g$ is reminiscent of the expansion of $P_{{\bf K}}g$ into a basis consisting of eigenfunctions for $\tau$. This also suggests that, perhaps, a  time-frequency approach to the similar open  questions concerning trilinear averages  will involve quadratic exponentials like $e^{i\lambda x^2}$, which are second order eigenfunctions for the rotations on the torus.   

Both theorems above will be consequences of the following very general harmonic analysis result, as explained in Section \ref{sec:transfer}. We will  use the notation
\begin{equation*}  
\label{e.dilate}
\text{Dil}_s^p h (x) = s^{-1/p}h(\frac{x}{s}) 
\end{equation*}
\begin{equation*}
\text{Mod}_\theta g(x)= e^{2\pi i \theta x} g(x). 
\end{equation*}

\begin{theorem} 
\label{t.osc1}  
Let $K:\R\to\R$ be an $L^2$ kernel satisfying the requirements:
\begin{equation}
\label{equation.1a11}
\widehat{K}\in C^{\infty}(\R\setminus\{0\}) 
\end{equation}
\begin{equation}
\label{equation.1a12}
|\widehat{K}(\xi)|\lesssim \min\{1,\frac{1}{|\xi|}\},\;\;\xi\not=0
\end{equation}
\begin{equation}
\label{equation.1a13}
|\frac{d^n}{d\xi^n}\widehat{K}(\xi)|\lesssim \frac{1}{|\xi|^n}\min\{|\xi|,\frac1{|\xi|}\},\;\;\xi\not=0,\;n\ge 1.
\end{equation}
Then for each $d=2^{1/n}$, $n\in \N$, each $f,g\in L^{\infty}(\R)$ with bounded support and  each finite sequence of integers $u_1<u_2<\ldots<u_J$ 
\begin{equation*}  
\label{e.returntimeo}
\|(\sum_{j=1}^{J-1}\sup_{k\in\Z\atop{u_j\le k<u_{j+1}}}|\int f(x+y)g(x-y)(\hbox{Dil}_{d^k}^1K({y})-\hbox{Dil}_{d^{u_{j+1}}}^1K({y}))dy|^2)^{1/2}\|_{L^{1,\infty}}
\end{equation*}
\begin{equation*}
\lesssim J^{1/4}\|f\|_{L^2}\|g\|_{L^2}
\end{equation*}
with the implicit constants depending only on $n$.
\end{theorem}

\begin{remark}
Due to the assumptions that $f$ and $g$ are bounded  and have bounded support it follows that for each $x\in\R$ the integral 
$$\int f(x+y)g(x-y)(\hbox{Dil}_{d^k}^1K({y})-\hbox{Dil}_{d^{u_{j+1}}}^1K({y}))dy$$ converges absolutely. The same remark applies to the following two theorems.
\end{remark}
\begin{remark}
The proof of this theorem is inspired by ideas from \cite{DTT}, \cite{DLTT} and  \cite{Lac1}. The same techniques can extend the theorem to a larger range for $p$ and $q$ and eliminate the dependence on $J$ of the bound in the above inequality. However, its current form suffices for our purposes. Moreover, the argument for a more general result as mentioned would have to be more technical and would require us to revisit most of the main results in \cite{DTT}, making the whole presentation much longer. One of the advantages of the current approach is that  it does  not rely on any type of interpolation, being a purely $L^2$ argument.

Our hope is that by keeping technicalities to a minimum, the whole argument will become easier to follow. The interested reader is certainly referred to \cite{DTT} for details on some of the results we are quoting here.
\end{remark}

As an immediate corollary of Theorem ~\ref{t.osc1} we get a  particular case of Lacey's inequality for the bilinear maximal function
\begin{corollary}[\cite{La2}] 
\label{corbilmax}
The following inequality holds for each $f,g\in L^2(\R)$
$$ \|\sup_{\epsilon>0}\frac{1}{\epsilon}\int_{|y|>\epsilon}|f(x+y)g(x-y)|dy\|_{1,\infty}\lesssim \|f\|_2\|g\|_2.$$
\end{corollary}
Theorem ~\ref{t.osc1} will follow from  two distinct results of dyadic analysis. The first one, Theorem ~\ref{t.oscpart}, is the particular case $d=2$ of the above and captures the main difficulty of the problem. The second one, Theorem ~\ref{t.square}, is a square function estimate and will be used to control error terms.

To understand better the {sec:transfer}connection between these three results we introduce some notation.
Let $x:(0,\infty)\to\C$. Let also $u_1<\ldots<u_{J}$ be as in Theorem ~\ref{t.osc1}  and define $a_1\le\ldots\le a_{J}$ such that $a_jn\le u_j<(a_j+1)n$. Then observe that 
\begin{align*}
(\sum_{j=1}^{J-1}\sup_{k\in\Z\atop{u_j\le k<u_{j+1}}}|x(\frac{k}{n})-x(\frac{u_{j+1}}{n})|^2)^{1/2}&\lesssim \sum_{i=0}^{n-1}(\sum_{j=1}^{J-1}\sup_{k\in\Z\atop{a_j\le k<a_{j+1}}}|x(k+\frac{i}{n})-x(a_{j+1}+\frac{i}{n})|^2)^{1/2}\\&+\sum_{i,j=0\atop{i\not=j}}^{n-1}(\sum_{k\in\Z}|x(k+\frac{i}{n})-x(k+\frac{j}{n})|^2)^{1/2}.
\end{align*} 

Using this inequality and a dilation argument, Theorem ~\ref{t.osc1} will follow immediately from the following two results.

\begin{theorem} 
\label{t.oscpart}  
Let $K:\R\to\R$ be an $L^2$ kernel satisfying ~\eqref{equation.1a11}, ~\eqref{equation.1a12} and ~\eqref{equation.1a13}. 
Then for each finite sequence of integers $u_1<u_2<\ldots<u_J$ and each $f,g\in L^{\infty}(\R)$ with bounded support
$$
\|(\sum_{j=1}^{J-1}\sup_{k\in\Z\atop{u_j\le k<u_{j+1}}}|\int f(x+y)g(x-y)(\hbox{Dil}_{2^k}^1K({y})-\hbox{Dil}_{2^{u_{j+1}}}^1K({y}))dy|^2)^{1/2}\|_{L^{1,\infty}}$$
$$
\lesssim J^{1/4}\|f\|_{L^2}\|g\|_{L^2},
$$
with some universal implicit constant.
\end{theorem}

\begin{theorem} 
\label{t.square}  
Let $K:\R\to\R$ be an $L^2$ kernel satisfying ~\eqref{equation.1a11}, ~\eqref{equation.1a12} and ~\eqref{equation.1a13} and the extra requirement
\begin{equation}
\label{equation.1a14}
|\widehat{K}(\xi)|\lesssim |\xi|,\;\;\xi\not=0.
\end{equation}
Then the  following inequality holds for each $f,g\in L^{\infty}(\R)$ with bounded support
\begin{equation*}  
\label{e.osc2}
\|(\sum_{k\in \Z}|\int f(x+y)g(x-y)\hbox{Dil}_{2^{k}}^1K({y})dy|^2)^{1/2}\|_{L^{1,\infty}}
\lesssim \|f\|_{L^2}\|g\|_{L^2}
\end{equation*}
with some universal  implicit constant.
\end{theorem}

In Section ~\ref{sec:transfer} we indicate how the result of Theorem ~\ref{t.osc1} can be transfered to a similar inequality in a dynamical system, and how this implies the convergence in theorems ~\ref{pontconvave} and ~\ref{pontconvser}. In Section ~\ref{sec:dis} we discretize the operator in Theorem ~\ref{t.oscpart} while the remaining sections are concerned with proving its boundedness. In the last section we briefly sketch how the same procedure can be applied to prove 
Theorem ~\ref{t.square}.

\section{Notation}
\label{sec:not}
In this section we set up some notation and terminology for the rest of the paper.
  
If $I$ is an interval then  $c(I)$ denotes the center of $I$, while $cI$ is the interval with the same center and length $c$ times the length of $I$.
By $ 1_{A}$ we denote the characteristic function of the set $A\subset \R$, while for any dyadic interval $I$, $$\chi_I(x)=\left(1+\frac{(x-c(I))^2}{|I|^2}\right)^{-1/2}.$$

A tile $P$ is a rectangle $P=I_P\times\omega_{P}$ such that  $I_P$ is a  dyadic interval and $|I_P|\cdot|\omega_P|=1$. A multitile $s$ is a box $s=I_s\times\omega_{s,1}\times\ldots\times\omega_{s,n}$ such that $I_s$ is a  dyadic interval, $|\omega_{s,1}|=\ldots|\omega_{s,n}|$ and  $|I_s|\cdot|\omega_{s,1}|=1$.

For each $E$ of finite measure,  $X(E)$ will denote  the set of all functions supported in $E$ with $\|f\|_{\infty}\le 1$, while $X_2(E)$ will denote  the $L^2$ normalized set of all functions supported in $E$ with $\|f\|_{\infty}\le |E|^{-1/2}$. Also $Mf(x)=\sup_{r>0}\frac{1}{2r}\int_{x-r}^{x+r} |f|(y)dy$ denotes the classical Hardy-Littlewood maximal function.

The notation $a\lesssim b$ means that $a\le cb$ for some universal constant $c$, while $a\sim b$ means that  $a\lesssim b$ and  $b\lesssim a$. Sometimes we will write $a\lesssim_{\hbox{parameters}} b$ to indicate that $a\le cb$ with $c$ depending only on the specified parameters.

\section{Pointwise convergence for averages and series}
\label{sec:transfer}
\subsection{Bounded oscillation implies convergence}

Assume we have a sequence $W_k$ of weighted operators defined on a dynamical system ${\bf X}=(X, \Sigma, m, \tau)$ by the formula
$$W_{k}(f,g)(x)=\sum_{n\in\Z}w_{k,n}f(\tau^{n}x)g(\tau^{-n}x),\;\;k\in\N.$$
\begin{lemma}
\label{oscimplconv}
If for some $f,g\in L^2(X)$
\begin{equation}
\label{nowcontrsihol}  
\|(\sum_{j=1}^{J-1}\sup_{k\in\Z\atop{u_j\le k<u_{j+1}}}|W_{k}(f,g)(x)-W_{u_{j+1}}(f,g)(x)|^2)^{1/2}\|_{L^{1,\infty}}
\lesssim J^{1/4}\|f\|_{L^2}\|g\|_{L^2}
\end{equation}
uniformly in $J$ and all finite  sequences of positive integers $u_1<u_2<\ldots<u_{J}$, then
$$\lim_{k\to\infty}W_{k}(f,g)(x)$$
exists for almost every $x\in X$.
\end{lemma}

\begin{proof}
To see this, assume for contradiction that the convergence does not hold. It follows that there is a measurable set $X'\subset X$ with $m(X')>0$ and some $\alpha>0$, such that for each $x\in X'$ 
$$\limsup_{k\to\infty}W_{k}(f,g)(x)-\liminf_{k\to\infty}W_{k}(f,g)(x)>\alpha.$$
 An elementary measure theoretic argument shows that one can then choose a subset $X''\subseteq X'$ of positive measure and a sequence of positive integers $(u_j)_{j\in\N}$ such that

$$\sup_{u_j\le k<u_{j+1}}|W_{k}(f,g)(x)-W_{u_{j+1}}(f,g)(x)|>\alpha,$$
for each $j\in \N$ and for each $x\in X''$. We immediately get that for each $J$

\begin{equation*}  
\|\left(\sum_{j=1}^{J}\sup_{u_j\le k<u_{j+1}}|W_{k}(f,g)(x)-W_{u_{j+1}}(f,g)(x)|^2\right)^{1/2}\|_{1,\infty}
  \ge J^{1/2}\alpha m(X''),
\end{equation*}
which contradicts inequality \eqref{nowcontrsihol}.
\end{proof}

\subsection{Proof of Theorem \ref{pontconvave}}

Let $M\in\N$ be  arbitrary. We apply Theorem \ref{t.osc1} to a $C^{\infty}$ kernel $K_{M}$ satisfying
$$1_{[0,1]}\le K_M\le 1_{[-\frac1M,1+\frac1M]}$$
Then we invoke standard transfer methods like in \cite{DTT} to get the following corollary.
\begin{corollary}
For each $d=2^{1/n}$, $n\in \N$,  each finite sequence of positive integers $u_1<u_2<\ldots<u_J$ and each $f,g\in L^2(X)$
\begin{equation*}  
\|(\sum_{j=1}^{J-1}\sup_{k\in\Z\atop{u_j\le k<u_{j+1}}}|(A_{d^k}(f,g)+E_{d^k,M}(f,g))-(A_{d^{u_{j+1}}}(f,g)+E_{d^{u_{j+1}},M}(f,g))|^2)^{1/2}\|_{{1,\infty}}
\end{equation*}
\begin{equation*}
\lesssim_{M,n} J^{1/4}\|f\|_{L^2}\|g\|_{L^2},
\end{equation*}
where for each real $r>1$ we denote
 $$A_r(f,g)(x)=\frac1{[r]}\sum_{0\le n\le r}f(\tau^nx)g(\tau^{-n}x)$$
while 
$$E_{r,M}(f,g)(x)=\sum_{n\in\Z}w_{r,n,M}f(\tau^nx)g(\tau^{-n}x)$$
is some error term with
\begin{equation}
\label{hcdhlah8`pdipoidk}
\sup_{r>1}\sum_{n\in\Z}|w_{r,n,M}|\lesssim M^{-1}.
\end{equation}
\end{corollary} 
By using this and Lemma \ref{oscimplconv} we find that there is a set $X_0$ of full measure such that for each $M\in\N$, each $d=2^{1/n}$ and each $x\in X_0$ the following limit exists
$$\lim_{k\to\infty}(A_{d^k}(f,g)(x)+E_{d^k,M}(f,g)(x)).$$

Fix now some $f,g\in L^{\infty}(X)$ with $\|f\|_{\infty},\|g\|_{\infty}\le 1$.
Based on \eqref{hcdhlah8`pdipoidk} we now get that for almost every $x\in X_0$
$$\limsup_{k\to\infty}A_{k}(f,g)(x)-\liminf_{k\to\infty}A_{k}(f,g)(x)\le$$ $$\le \limsup_{k\to\infty}(A_{d^k}(f,g)(x)+E_{d^k,M}(f,g)(x))-\liminf_{k\to\infty}(A_{d^k}(f,g)(x)+E_{d^k,M}(f,g)(x))+$$
$$+\limsup_{k\to\infty}E_{k}(f,g)(x)-\liminf_{k\to\infty}E_{k}(f,g)(x)+10n^{-1}$$
$$\lesssim n^{-1}+M^{-1}.$$
Since we can take  $M$ and $n$ to be as large as we want, it follows that 
$$\lim_{k\to\infty}A_{k}(f,g)(x)$$
exits almost everywhere.

\subsection{Proof of Theorem \ref{pontconvser}}

Since the Hilbert kernel involved in Theorem \ref{pontconvser} is not integrable, the route towards proving convergence in this case poses some extra difficulties. As a consequence, we will present a more detailed  argument in this case.

The point is again to prove the almost everywhere convergence of the series
$$\sideset{}{'}\sum_{n=-N}^{N}\frac{f(\tau^nx)g(\tau^{-n}x)}{n}$$
along lacunary sequences and then invoke the boundedness of both $f$ and $g$ to get the convergence along the full sequence of positive integers. For simplicity we choose to present the argument in the particular case $N=2^k$.

Let $M\in\N$ be  arbitrary. We apply Theorem \ref{t.osc1} to a $C^{\infty}$ kernel $K_{M}$ satisfying
$$K_M(x)=\frac{1}{x}\;\;\hbox{for\;}|x|\ge 1$$
$$|K_M1_{[-1,1]}|\le 2\times 1_{1-\frac1M\le |x|\le 1}.$$
Introduce the discrete kernels $H_{k,M}:\R\to\R,\; k\ge 1$,  defined by the following
$$H_{k,M}(x)=\sum_{-2^k\le i\le 2^{k}-1}1_{[i,i+1)}(x)\frac{1}{2^k}K_{M}(\frac{i}{2^k})+\sum_{i\in\Z\setminus [-2^k,2^k-1]}1_{[i,i+1)}(x)\frac{1}{i}.$$

Let $k_1<k_2<\ldots<k_J$  be an arbitrary sequence of positive integers. The fact that the terms of the sequence are positive is a crucial fact, exploited in the following. Indeed, we note that
$$|H_{k,M}(y)-\hbox{Dil}_{2^k}^1K_M(y)|\lesssim_{M} \begin{cases}&\frac1{2^{2k}},\;\; |y|\le 2^k\\&\frac{1}{y^2},\;\;|y|\ge 2^k\end{cases}.$$ 
From the boundedness of the maximal averages in Corollary \ref{corbilmax} we deduce that
\begin{align*}  
\|(\sum_{k\ge 1}|\int &f(x+y)g(x-y)(H_{k,M}({y})-\hbox{Dil}_{2^k}^1K_M(y))dy|^2)^{1/2}\|_{1,\infty}\\&
\lesssim_{M} \|\sup_{\epsilon>0}\frac{1}{\epsilon}\int_{|y|>\epsilon}|f(x+y)g(x-y)|dy\|_{1,\infty}\\&\lesssim \|f\|_2\|g\|_2.
\end{align*}
As a consequence of this and Theorem \ref{t.osc1} applied to $K_M$, we get that 
\begin{equation*}  
\|(\sum_{j=1}^{J-1}\sup_{k_j\le k<k_{j+1}}|\int f(x+y)g(x-y)(H_{k,M}({y})-H_{k_{j+1},M}({y}))dy|^2)^{1/2}\|_{1,\infty}
\end{equation*}
\begin{equation}  
\label{e.returntime14o}
\lesssim_{M} J^{1/4}\|f\|_2\|g\|_2.
\end{equation}

The next step consists of transferring  ~\eqref{e.returntime14o} to the integers. By considering functions  like $f:\R\to\R$ with $$f(x)=\begin{cases}\phi([x])&:\quad [x]+\frac14\le x\le [x]+\frac12 \\ \hfill  0&:\quad \text{otherwise}\end{cases}$$ and  $g:\R\to\R$ with $$g(x)=\begin{cases}\psi([x])&:\quad [x]+\frac14\le x\le [x]+\frac12 \\ \hfill  0&:\quad \text{otherwise}\end{cases},$$ we get that for each $\phi,\psi:\Z\to\Z$ with finite support 
$$
\|(\sum_{j=1}^{J-1}\sup_{k_j\le k<k_{j+1}}|\sum_{b\in\Z} \phi(a+b)\psi(a-b)(H_{k,M}(b)-H_{k_{j+1},M}(b))|^2)^{1/2}\|_{l^{1,\infty}(\Z)}
$$
\begin{equation}
\label{e.returntime16o}\lesssim_{M} J^{1/4}\|\phi\|_{l^2(\Z)}\|\psi\|_{l^2(\Z)}.
\end{equation}

For each $k\ge 1$ introduce the  kernels $A_{k,M}:\Z\to\Z$ and $S_{k,M}:\Z\to\Z$ defined by
$$A_{k,M}(i)=\begin{cases}&H_{k,M}(i),\;\;-2^k\le i\le 2^k\\&0,\;\;\hbox{otherwise}\end{cases},$$

$$S_{k,M}(i)=\begin{cases}&\frac1i,\;\;-2^k\le i\le 2^k,\;i\not=0\\&0,\;\;\hbox{otherwise}\end{cases},$$
and note that for each $k<k'$
$$H_{k,M}-H_{k',M}=O_{k,M}-O_{k',M}:=(A_{k,M}-S_{k,M})-(A_{k',M}-S_{k',M}).$$ Thus ~\eqref{e.returntime16o} gives 
\begin{equation*}  
\|(\sum_{j=1}^{J-1}\sup_{k_j\le k<k_{j+1}}|\sum_{b\in\Z} \phi(a+b)\psi(a-b)(O_{k,M}(b)-O_{k_{j+1},M}(b))|^2)^{1/2}\|_{l^{1,\infty}(\Z)}
\lesssim_{M} J^{1/4}\|\phi\|_{l^2(\Z)}\|\psi\|_{l^2(\Z)}.
\end{equation*}

Standard transfer to a dynamical system ${\bf X}=(X,\Sigma, \mu,\tau)$, as in \cite{DTT}, leads to 
\begin{equation*}  
\|(\sum_{j=1}^{J-1}\sup_{k_j\le k<k_{j+1}}|\sum_{n\in\Z} f(\tau^nx)g(\tau^{-n}x)(O_{k,M}(n)-O_{k_{j+1},M}(n))|^2)^{1/2}\|_{1,\infty}
\end{equation*}
\begin{equation}
\label{e.returntime18o}
\lesssim_{M} J^{1/4}\|f\|_2\|g\|_2.
\end{equation}

By invoking Lemma \ref{oscimplconv} it follows that if $\|f\|_{L^{\infty}},\|g\|_{L^{\infty}}\le 1$
$$\lim_{k\to\infty}\sum_{n\in\Z} f(\tau^nx)g(\tau^{-n}x)O_{k,M}(n)$$
exists for almost every $x\in X$.
Finally
$$\limsup_{k\to\infty}\sum_{n\in\Z} f(\tau^nx)g(\tau^{-n}x)S_{k,M}(n)-\liminf_{k\to\infty}\sum_{n\in\Z} f(\tau^nx)g(\tau^{-n}x)S_{k,M}(n)\le $$
$$\le \limsup_{k\to\infty}\sum_{n\in\Z} f(\tau^nx)g(\tau^{-n}x)O_{k,M}(n)-\liminf_{k\to\infty}\sum_{n\in\Z} f(\tau^nx)g(\tau^{-n}x)O_{k,M}(n)+$$
$$+\limsup_{k\to\infty}\sum_{n\in\Z} f(\tau^nx)g(\tau^{-n}x)A_{k,M}(n)-\liminf_{k\to\infty}\sum_{n\in\Z} f(\tau^nx)g(\tau^{-n}x)A_{k,M}(n)\lesssim$$
$$\lesssim CM^{-1}.$$
Since this holds for arbitrary $M$, we get that
$$\lim_{k\to\infty}\sideset{}{'}\sum_{n=-2^k}^{2^k}\frac{f(\tau^nx)g(\tau^{-n}x)}{n}$$
exists almost everywhere.

\section{Discretization}
\label{sec:dis}
\begin{definition}
A set $\G'$ of (not necessarily dyadic) intervals is called a grid if
\begin{itemize}
\item   $I,I'\in \G'$ and $|I'|\le |I|$ imply that  either $I'\subset I$ or $I\cap I'=\emptyset$
\item  $I\in\G'$ implies that $|I|=2^k$, for some $k\in\Z$.
\end{itemize}
\end{definition}

The standard dyadic grid is
\begin{equation}
\label{ctstangrid}
\B=\{\;[2^il,2^i(l+1)]:\;i,l\in\Z\}.
\end{equation}
We will also be interested in a  more general type of grid. For each odd integer $N\ge 3$, each $0\le t\le N-2$ and $0\le L\le N-1$, the collection
$$\G_{N,t,L}:=\left\{\left[2^{i}\left(l+\frac{L}{N}\right),2^{i}\left(l+\frac{L}{N}+1\right)\right]\;:i\equiv t\pmod {N-1},\;l\in\Z\right\}$$
is a  grid, as it easily follows from the fact that $2^{N-1}\equiv 1\pmod N.$ 
We note that for each fixed $N$ the grids $\G_{N,t,L}$ are pairwise disjoint, for $0\le t\le N-2$ and $0\le L\le N-1$.

For each $L^2$ kernel $K$, each $f,g\in L^{\infty}(\R)$ with bounded support and each sequence ${\bf U}=(u_j)_{j=1}^{J}$ define the quantity

$$O_{K,{\bf U}}(f,g)(x)$$
$$=(\sum_{j=1}^{J-1}\sup_{k\in\Z\atop{u_j\le k<u_{j+1}}}|\int f(x+y)g(x-y)(\hbox{Dil}_{2^k}^1K({y})-\hbox{Dil}_{2^{u_{j+1}}}^1K({y}))dy|^2)^{1/2}.$$

In the following we indicate how to discretize it.
Choose  $\eta:\R\to\R$ such that $\widehat{\eta}$ is a $C^{\infty}(\R\setminus\{0\})$ function which equals $\lim_{\xi\to 0^{+}}\widehat{K}(\xi)$ on $\left(0,\frac12\right]$, $\lim_{\xi\to 0^{-}}\widehat{K}(\xi)$ on $\left[-\frac12,0\right)$ and $0$ outside  $\left[-1,1\right]$. The two limits exist due to the fact that $|\frac{d}{d\xi}K(\xi)|\lesssim 1$ for $\xi\not=0$.
 It suffices to prove
\begin{equation}
\label{e.osc5}
\|\|O_{\eta,{\bf U}}(f,g)\|_{1,\infty}\lesssim J^{1/4}\|f\|_{2}\|g\|_{2}
\end{equation}
\begin{equation}
\label{e.osc6}
\|\|O_{K-\eta,{\bf U}}(f,g)\|_{1,\infty}\lesssim J^{1/4}\|f\|_{2}\|g\|_{2}.
\end{equation}
The proofs for the above inequalities will follow from a more general principle, as explained below. The crucial property of the multiplier $\widehat{K-\eta}$ that will be used later is the following 
\begin{equation}
\label{e:osc7}
|\frac{d^n}{d{\xi}^n}\widehat{K-\eta}(\xi)|\lesssim \frac{1}{|\xi|^n}\min\{|\xi|,\frac{1}{|\xi|}\},\;\;n\ge 0.
\end{equation}
Note that the additional inequality $|\widehat{K-\eta}(\xi)|\lesssim |\xi|$ for $\xi\not=0$ is  a consequence of the fact that $|\frac{d}{d\xi}\widehat{K}(\xi)|\lesssim 1$ for $\xi\not=0$.
Write 
\begin{equation}
\label{e:osc8}
\widehat{K-\eta}(\xi)=\sum_{j=-\infty}^{\infty}\widehat{K-\eta}(\xi)q(\frac{\xi}{2^j}),
\end{equation}
 where $q$ is some Schwartz function supported in the annulus $\frac12<|\xi|<2$ such that 
$$\sum_{j\in \Z}q(\frac{\xi}{2^j})=1,\;\;\xi\not=0.$$ Define $g_j=\widehat{K-\eta}(\xi)q(\frac{\xi}{2^j})$. As a consequence of ~\eqref{e:osc7} we get that both the function $\eta$ and  each function $\operatorname{Dil}_{2^j}^1\check{g_j}$ will satisfy 
$$|{\eta}(x)|\lesssim_{M} \frac{1}{(1+|x|)^{M}}$$
$$|\frac{1}{2^j}\check{g_j}(\frac{x}{2^j})|\lesssim_{M} \frac{2^{-|j|}}{(1+|x|)^{M}},$$ 
for all $M\ge 0$ and $x\in\R$, uniformly in $j\in\Z$.
Moreover, each of the above functions has the Fourier transform  constant on both $\{0<\xi\le \frac12\}$ and $\{-\frac12\le \xi< 0\}$, as well as on $\{|\xi|\ge 2\}$. 

For each $j\in \Z$ define the shifted sequence $\U'_j=\U-j$. Since the operators $O_{\check{g_j},\U}$ and $O_{\operatorname{Dil}_{2^{j}}^{1}\check{g_j},\U'_j}$ coincide, inequalities   ~\eqref{e.osc5} and  ~\eqref{e.osc6} will immediately follow if we prove that
\begin{equation}
\label{en.2}
\|O_{\check{\theta},\U}(f,g)\|_{1,\infty}\lesssim J^{1/4}\|f\|_{2}\|g\|_{2},
\end{equation}
uniformly in all $C^{\infty}(\R\setminus\{0\})$ functions $\theta$ which are constant on both $\{0<\xi\le \frac12\}$ and $\{-\frac12\le \xi< 0\}$, as well as on $\{|\xi|\ge 2\}$, and which satisfy
\begin{equation}
\label{gdhgfhgtr}
|\check{\theta}(x)|\lesssim_{M} \frac1{(1+|x|)^M}
\end{equation}
for all $M\ge 0$ and $x\in\R$.

By a similar reduction we can assume instead that  $\theta$ is constant on $\{|x|\le 1000\}$ and on $\{|x|\ge 4000\}$. This extra assumption will serve later for the purpose of creating disjointness of some sort between multitiles. Note that 
$$\theta(2^k\xi)=\sum_{i\ge k}\theta_i(\xi)
$$ 
where $\theta_i(\xi)=\theta(2^{i}\xi)-\theta(2^{i+1}\xi)$ is supported in the annulus $\{500\times2^{-i}\le |\xi|\le 4000\times2^{-i}\}$.

Pick a Schwartz function $\psi$ such that $\widehat{\psi}$ is supported in $[0,\frac25]$ and satisfies the following property for every $\xi\in\R$
$$\sum_{l\in \Z}\left|\widehat{\psi}\left(\xi-\frac{l}5\right)\right|^2=1.
$$
For each scale $i$ use the following expansion for both $f$ and $g$, valid in every $L^p$ norm, $1<p<\infty$ 

$$f=\sum_{m,l\in \Z}\langle f,\psi_{i,m,\frac{l}5}\rangle\psi_{i,m,\frac{l}5},$$$$g=\sum_{m,l\in \Z}\langle g,\psi_{i,m,\frac{l}5}\rangle\psi_{i,m,\frac{l}5},$$where 
$$\psi_{i,m,l}(x)=2^{-\frac{i}{2}}\psi(2^{-i}x-m)e^{2\pi i2^{-i}xl}.$$ The fundamental properties of $\psi_{i,m,l}$ that will be used in the following are
\begin{equation}
\label{osc15}
|\frac{d^n}{d x^n}\operatorname{Mod}_{-l2^{-i}}\psi_{i,m,l}(x)|\lesssim_{M}2^{(-1/2-n)i}\frac1{(1+|2^{-i}x-m|)^M},\;\;M\ge 0
\end{equation}
\begin{equation}
\label{osc16}
\operatorname{supp}\widehat{\psi}_{i,m,l}\subset [\frac{l}{5}2^{-i},(\frac{l}{5}+1)2^{-i}].
\end{equation}
The above properties can be summarized by saying that the tile $[m2^i,(m+1)2^i]\times [\frac{l}{5}2^{-i},(\frac{l}{5}+1)2^{-i}]$ is a Heisenberg box of $\psi_{i,m,l}$.

With the notation 
$$\varphi_{i,\vec{m},\vec{l}}(x):=2^{i/2}\int\psi_{i,m_1,\frac{l_1}5}(x+y)\psi_{i,m_2,\frac{l_2}5}(x-y)\check{\theta_i}(y)dy$$
it follows that
$$O_{\check{\theta},\U}(f,g)(x)=$$
$$(\sum_{j=1}^{J-1}\sup_{k\in\Z\atop{k_j\le k<k_{j+1}}}|\sum_{i=k}^{k_{j+1}-1}\sum_{\vec{m},\vec{l}\in\Z^2}2^{-i/2}\langle f,\psi_{i,m_1,\frac{l_1}5}\rangle\langle g,\psi_{i,m_2,\frac{l_2}5}\rangle\varphi_{i,\vec{m},\vec{l}}(x)|^2)^{1/2}$$
Triangle's inequality then shows that Theorem ~\ref{t.oscpart} follows once  we prove that
$$
\|\sum_{m\ge 0}(\sum_{j=1}^{J-1}\sup_{k\in\Z\atop{u_j\le k<u_{j+1}}}|\sum_{i=k}^{u_{j+1}-1}\sum_{\vec{m},\vec{l}\in\Z^2\atop{|m_1-m_2|=m}}2^{-i/2}\langle f,\psi_{i,m_1,\frac{l_1}5}\rangle\langle g,\psi_{i,m_2,\frac{l_2}5}\rangle\varphi_{i,\vec{m},\vec{l}}|^2)^{1/2}\|_{1,\infty}$$
\begin{equation}
\label{osc13}
\lesssim J^{1/4}\|f\|_2\|g\|_2.
\end{equation}

The computations that follow are meant to reveal the decay and localization of the functions $\varphi_{i,\vec{m},\vec{l}}$.

We first observe that since 
$$\widehat{\varphi}_{i,\vec{m},\vec{l}}(\xi)=\int\widehat{\psi}_{i,m_1,\frac{l_1}5}(\eta)\widehat{\psi}_{i,m_2,\frac{l_2}5}(\xi-\eta)\theta(2\eta-\xi)d\eta,$$
it turns out that
\begin{equation}
\label{e.osc10}
\operatorname{supp}\widehat{\varphi}_{i,\vec{m},\vec{l}}\subseteq \operatorname{supp} 
\widehat{\psi}_{i,m_1,l_1}+\operatorname{supp}\widehat{\psi}_{i,m_2,l_2}\subset [\frac{l_3}{5}2^{-i},(\frac{l_3}{5}+1)2^{-i}],
\end{equation}
where from now on we will denote 
\begin{equation}
\label{e.osc19}
l_3=l_1+l_2.
\end{equation}

We next observe that for each $M\ge 0$
\begin{align*}
|2^{i/2}\varphi_{i,\vec{m},\vec{l}}(2^ix)|&\lesssim \int |\psi(x+y-m_1)\psi(x-y-m_2)\check{\theta}(y)|dy\\&\lesssim_{M} \frac{1}{(1+|x-m_1|)^M(1+|x-m_2|)^M}\\&\lesssim_{M} \frac1{(1+|m_1-m_2|)^{M}(1+|x-\frac{m_1+m_2}{2}|)^M}.
\end{align*}
A similar estimate holds for all derivatives, and we conclude that
\begin{equation}
\label{e.decay1}
|\frac{d^n}{d x^n}\operatorname{Mod}_{-\frac{l_3}{5}2^{-i}}\varphi_{i,\vec{m},\vec{l}}(x)|\lesssim_{M} 2^{(-1/2-n)i}\frac1{(1+|m_1-m_2|)^{M}(1+|2^{-i}x-\frac{m_1+m_2}{2}|)^M},
\end{equation}
for each $n,M\ge 0.$

We thus see that $\varphi_{i,\vec{m},\vec{l}}$ satisfies the same type of properties as $\psi_{i,m,l}$, with some extra uniform decay in $m=|m_1-m_2|.$ In particular, if $I_{i,\vec{m}}$ is one of the (at most two) dyadic intervals of length $2^i$ which contains $\frac{m_1+m_2}{2}$, then $I_{i,\vec{m}}\times[\frac{l_3}{5}2^{-i},(\frac{l_3}{5}+1)2^{-i}]$ is certainly a Heisenberg box for $\varphi_{i,\vec{m},\vec{l}}.$ 

Finally, since
$$\varphi_{i,\vec{m},\vec{l}}(x)=\int\widehat{\psi}_{i,m_1,\frac{l_1}5}(\xi_1)\widehat{\psi}_{i,m_2,\frac{l_2}5}(\xi_2)\theta_i(\xi_2-\xi_1)e^{2\pi ix(\xi_1+\xi_2)}d\xi_1 d\xi_2,$$
it follows that in order for $\varphi_{i,\vec{m},\vec{l}}$ to be non identically equal to 0 we must have
\begin{equation}
\label{e.osc9}
10^2< |\frac{l_1}{5}-\frac{l_2}{5}|< 10^{4}.
\end{equation}

The next reduction concerns the fact that there is a finite universal  set $E\in \Z$, such that every every $l_1,l_2,l_3$ satisfying ~\eqref{e.osc19}  and ~\eqref{e.osc9} will also satisfy the following 
\begin{equation}
\label{pat1}
l_2=l_1+e
\end{equation}
for some $e\in E.$
We will restrict the summation in ~\eqref{osc13} to those vectors  $\vec{l}$ satisfying ~\eqref{pat1} for some fixed $e$.

For each $\vec{l}\in\Z^2$ and each $i\in\Z$ define the cubes $Q_{\vec{l},i}=\prod_{j=1}^{3}[\frac{l_j}{5}2^{-i},(\frac{l_j}5+1)2^{-i}]$.
 Note  that  every  frequency interval $[2^{-i}\frac{l_j}5,2^{-i}(\frac{l_j}5+1)]$ of a cube $Q_{\vec{l},i}$ belongs to one of the grids  $\G_{5,t,L}$, $0\le t\le 3$, $0\le L\le 4.$
This allows us to further restrict the summation in ~\eqref{osc13} to those $\vec{l}$ and $i$ for which each 
interval $[2^{-i}\frac{l_j}5,2^{-i}(\frac{l_j}5+1)]$ is in a fixed grid depending on $j$.

Denote by $\D'$  the union over all $i\in\Z$ and all $l_1,l_2,l_3\in \Z$ which are subject to all the indicated restrictions, of the set of all the cubes $Q_{\vec{l},i}$. For each $m\ge 0$ introduce the set of  generalized multitiles $\P_{m}$ to be
$$\{P=\prod_{j=1}^{3}I_{P,j}\times Q_{\vec{l},i}:\;I_{P,j}=[m_j2^i,(m_j+1)2^i],\;j\le 2,\; I_{P,3}=I_{i,\vec{m}},\;Q_{\vec{l},i}\in\D',\;|m_1-m_2|=m\}.
$$ Each such multitile can be thought of as the product of $3$ tiles $P_j:=I_{P,j}\times [l_j2^{-i},(l_j+1)2^{-i}]$. For each multitile $P$ as above define
$\psi_{P,j}=\psi_{i,m_j,l_j}$ for $j=1,2$ and  $\psi_{P,3}=\varphi_{i,\vec{m},\vec{l}}$. Thus ~\eqref{osc13} will follow if we prove that
\begin{equation*}
\|(\sum_{j=1}^{J-1}\sup_{k\in\Z\atop{u_j\le k<u_{j+1}}}|\sum_{P\in\P_m\atop{2^k\le |I_P|<2^{u_{j+1}-1}}}|I_P|^{-1/2}\langle f,\psi_{P,1}\rangle\langle g,\psi_{P,2}\rangle\psi_{P,3}|^2)^{1/2}\|_{1,\infty}
\end{equation*}
\begin{equation*}
\lesssim (1+m)^{-2}J^{1/4}\|f\|_2\|g\|_2,
\end{equation*}
where $|I_P|$ is defined as being the common value of all $|I_{P,j}|$.

As it will easily follow from our later analysis, it is enough to prove the above for $m=0$. The extra decay in $m$ for the other terms will be a consequence of the extra decay in ~\eqref{e.decay1}. Once we restrict attention to the case $m=0$, we remark that $I_{P,1}=I_{P,2}=I_{P,3}$, however this observation will not change or simplify our argument.

A last harmless reduction consists of  sparsifying the set of time dyadic intervals. More precisely, we will assume that if $\frac{|I'|}{|I|}>1$ then $\frac{|I'|}{|I|}\ge 2^{\Delta}$, where $\Delta$ is a sufficiently large constant to be chosen later  ($\Delta=1000$ will certainly suffice). 

We now summarize all the various reductions made so far in the following theorem, which implies Theorem ~\ref{t.oscpart}.

\begin{theorem}
\label{t.oscprincmodel}
Let $\G$, $\G_1$, $\G_2$ and $\G_3$ be four grids with $\G$ satisfying
\begin{equation}
\label{e.thmgrid1}
\G\subset\B
\end{equation}
\begin{equation}
\label{e.thmgrid2}
I,I'\in\G\Rightarrow\max\{|I||I'|^{-1},|I'||I|^{-1}\}\ge 2^{\Delta}. 
\end{equation}
Let $e$ be a  number with $10^2\le |e|\le 10^5$ and define
$$\D=\{\prod_{j=1}^{3}[\frac{l_j}{5}2^i,(\frac{l_j}{5}+1)2^i]\in\prod_{j=1}^{3}\G_j:\; l_2=l_1+e,\;l_3=l_1+l_2\}.$$
Define also the  set of multitiles 
$$\S=\{s=I_{s}\times Q_{s}: I_s\in\G,\:\: Q_s\in\D\;\text{with sidelength}\;\frac{1}{|I_s|}\}.$$
Assume that  each multitile $s=I_s\times\prod_{j=1}^{3}\omega_{s,3}\in\S$ is associated with three functions $(\psi_{s,j})_{j=1}^{3}$ satisfying
\begin{equation}
\label{e.osc23}
|\frac{d^n}{d x^n}\operatorname{Mod}_{-c(\omega_{s,j})}\psi_{s,j}(x)|\lesssim_{n,M} |I_s|^{-1/2-n}\chi_{I_s}^M,\;\;n,M\ge 0
\end{equation}
\begin{equation}
\label{e.osc24}
\operatorname{supp}\widehat{\psi}_{s,j}\subset \omega_{s,j},
\end{equation}
for each $j=1,2,3.$ 

Then  for each $f,g\in L^2(\R)$ and each finite sequence of integers $\U:=u_1<u_2<\ldots<u_{J}$ we have the estimate
\begin{equation}
\label{e.form1}
\|(\sum_{j=1}^{J-1}\sup_{k\in\Z\atop{u_j\le k<u_{j+1}}}|\sum_{s\in\S\atop{2^k\le |I_s|<2^{{u_{j+1}}}}}|I_s|^{-1/2}\langle f,\psi_{s,1}\rangle\langle g,\psi_{s,2}\rangle\psi_{s,3}|^2)^{1/2}\|_{1,\infty}\lesssim J^{1/4}\|f\|_2\|g\|_2,
\end{equation}
with the implicit constant depending only on the implicit constants in ~\eqref{e.osc23}.
\end{theorem}
 
For each $1\le j\le J-1$ let $\kappa_j:\R\to \{u_j,u_j+1,\ldots,u_{j+1}-1\}$ be some arbitrary stopping time. For each $s\in\S$ with $2^{u_1}\le |I_s|<2^{u_{J}}$ define 
$$\phi_{s,j}(x)=\psi_{s,j}(x)$$ for $j\in\{1,2\}$ and $$\phi_{s,3}(x)=\psi_{s,3}1_{2^{\kappa_j(x)}\le |I_s|<2^{u_{j+1}}}(x),$$ where $j$ is the unique integer such that $2^{u_j}\le|I_s|<2^{u_{j+1}}$. An equivalent formulation for ~\eqref{e.form1} that we will sometimes find easier to handle  is
\begin{equation}
\label{e.form2}
\|(\sum_{j=1}^{J-1}|\sum_{s\in\S\atop{2^{u_j}\le|I_s|<2^{u_{j+1}}}}|I_s|^{-1/2}\langle f,\phi_{s,1}\rangle\langle g,\phi_{s,2}\rangle\phi_{s,3}|^2)^{1/2}\|_{1,\infty}\lesssim J^{1/4}\|f\|_2\|g\|_2.
\end{equation}

From now on we will fix the collection $\S$, the wave packets $\psi_{s,j}$, the sequence $\U$ and the stopping times $\kappa_j$, and we will implicitly assume that for each $s\in\S$ we have $2^{u_1}\le |I_s|<2^{u_{J}}$.

\section{The combinatorics of the multitiles}
In this section we start by  defining a relation of order between multitiles. This will alow us to split $\S$ into structured collections like trees and forests. The model sum restricted to each tree is essentially a Littlewood-Paley dyadic decomposition modulated with  a frequency from the frequency interval of the top of the tree. Estimates for the model sum restricted to each such tree involves classical Calderon-Zygmund theory. This will be seen in Section ~\ref{sec:single}. The modern time-frequency side of the whole approach manifests in the fact that $\S$ consists of many trees, modulated with possibly different frequencies. The goal of this section is to prove that the model sums corresponding to distinct trees are almost orthogonal, a principle quntified in various Bessel type inequalities. The almost orthogonality will follow if the trees are selected to be strongly disjoint, a combinatorial property which, as the name suggests, is stronger than mere disjointness.
  
\begin{definition}
 For two multitiles $s,s'\in\S$ we write 
\begin{itemize}
\item $s_j'<s_j$ if $I_{s'}\varsubsetneq I_s$ and $\omega_{s,j}\varsubsetneq\omega_{s',j}$
\item $s_j'\le s_j$ if  $s_j'<s_j$ or $s_j=s_j'$
\item $s_j'\lesssim s_j$ if $I_{s'}\subseteq I_s$ and $\omega_{s,j}\subset 10e\omega_{s',j}$
\item  $s_j'\lesssim' s_j$ if    $s_j'\lesssim s_j$ and  $10\omega_{s,j}\cap 10\omega_{s',j}=\emptyset$, 
\end{itemize}
where $e$ is the constant in Theorem ~\ref{t.oscprincmodel}.
\end{definition}

\begin{lemma}
\label{lem:7}
Given any two multitiles $s,s'\in\S$ such that $s_i'<s_i$ for some $i\in\{1,2,3\}$, it follows that $s_j'\lesssim' s_j$ for each $j\in\{1,2,3\}\setminus\{i\}.$ 
\end{lemma}
\begin{proof}
We illustrate the proof on a particular case which is otherwise representative for the general argument.
Consider some $Q_s=\prod_{j=1}^3[\frac{l_j}{5}2^i,(\frac{l_j}{5}+1)2^i]$, $Q_{s'}=\prod_{j=1}^3[\frac{l_j'}{5}2^i,(\frac{l_j'}{5}+1)2^i]$ and assume that $s_1'<s_1$. This implies that 
\begin{equation}
\label{e.rank1}
|l_1'-2^{i-i'}l_1|\le 5.
\end{equation}
 
Assume now for contradiction  that $\omega_{s,2}\nsubseteq 10e\omega_{s',2}$. This in turn implies that 
\begin{equation}
\label{e.rank2}
|l_2'-2^{i-i'}l_2|> 10e.
\end{equation}
Inequalities ~\eqref{e.rank1} and ~\eqref{e.rank2} together with the fact that $l_2-l_1=l_2'-l_1'=e$ will immediately lead to the contradiction.

Assume next for contradiction  that $10\omega_{s,2}\cap 10\omega_{s',2}\not=\emptyset$. This in turn implies that 
\begin{equation}
\label{e.rank3}
|l_2'-2^{i-i'}l_2|\le 50. 
\end{equation}
Inequalities ~\eqref{e.rank1} and ~\eqref{e.rank3} together with the fact that  $l_2-l_1=l_2'-l_1'=e$ 
 will immediately lead to the contradiction.
\end{proof}

Let $i,j\in\{1,2,3\}$.
\begin{definition}
A collection of multitiles $\T\subset \S$ is called an $i$-tree with top $T\in \S$ if 
$$s_i<T_i\;\;\text{for each}\; s\in \T\setminus\{T\}.$$
We say that $\T$ is a tree if it is i-tree for some i. 
We call the tree $\T$ j-lacunary if 
$$s_j\lesssim' T_j\;\;\text{for each}\; s\in \T\setminus \{T\}.$$
\end{definition}
\begin{remark}
A tree does not necessarily contain its top.
By Lemma ~\ref{lem:7} we know that a tree is $j$-lacunary if and only if it is an $i$-tree for some $i\not=j$.  Note also that each tree can contain at most one multitile  with a given time interval.
\end{remark}

A very useful tool for proving estimates  for a single tree $\T$ is its size, a quantity which encodes the BMO properties of the model sum associated with $\T$. 
\begin{definition}
Consider some $j\in\{1,2,3\}$ and a finite subset of multitiles $\S'\subset\S$. Assume that each $s\in\S'$ is associated with a function $F_s:\R\to\C$. Define the  $j$-size of $\S'$ relative to the collection $(F_s)$ by the formula 
$$\size_j(\S')=\sup_{\T}\left(\frac{1}{|I_T|}\sum_{s\in \T}|\langle F_s, \phi_{s,j}\rangle|^2\right)^{\frac{1}{2}}$$ where the supremum is taken over all the trees $\T \in \S'$ of lacunary type $j$. 
\end{definition} 
To simplify notation, we will not add $(F_s)$ as a superscript of $\size$, since $(F_s)$ will always be clear from the context. Actually the 1-size and the 2-size will always be understood with respect to the functions $f$ and $g$, respectively. 

Note  that  size is a monotone function with respect to $\S'$.
The following lemma shows how to estimate tree paraproducts by using the size.
\begin{lemma}
\label{l:328}
If $\T$ is a tree then
$$\sum_{s\in \T}|I_s|^{-1/2}\prod_{i=1}^{3}|\langle F_s, \phi_{s,i}\rangle|\lesssim |I_T|\prod_{i=1}^{3}\size_i(\T).$$
\end{lemma}
\begin{proof}
Assume $\T$ is an $i$-tree. Apply $l_2$ estimates  for the terms corresponding to $j,k\in\{1,2,3\}\setminus \{i\}$ and $l_{\infty}$ for the terms corresponding to $i$.
\end{proof}
\begin{definition}
Let $j\in\{1,2,3\}$. Two trees $\T$ and $\T'$ with tops $T$ and $T'$ are said to be strongly $j$-disjoint if 
\begin{itemize}
\item
$s_j\cap s_j'=\emptyset$ for each $s\in \T\cup\{T\}, s'\in \T'\cup\{T'\}$ 

\item
whenever $s\in \T\cup\{T\}, s'\in \T'\cup\{T'\}$ are such that $\omega_{s,j}\subsetneq \omega_{s',j}$, then one has $I_T\cap I_{s'}=\emptyset$, and similarly with $\T$ and $\T'$ reversed.
\end{itemize}
A collection of trees is called mutually strongly $j$-disjoint if each two trees in the collection are strongly $j$-disjoint and each tree is $j$-lacunary.
\end{definition}
\begin{remark}
\label{rem:md123}
Each subset $\S_{\T}\subset\T$ of an $i$-tree $\T$ can be decomposed in a unique way as the disjoint union  of $i$-trees $\T'$ containing their tops 
$$\S_{\T}=\bigcup_{\T'\in D(\S_{\T})}\T',$$
such that these tops have  pairwise disjoint time intervals. Precisely, the tops of the trees in the collection $D(\S_{\T})$ will be the maximal multitiles in $\S_{\T}$ with respect to the order $<$. If one has a collection $\F$ of  mutually strongly $j$-disjoint trees $\T$  and each $\T$ is associated with a subset  $\S_{\T}\subset\T$, then the decomposition
$$\F=\bigcup_{\T\in\F}\bigcup_{\T'\in D(\S_{\T})}$$
 gives rise to a family of  mutually strongly $j$-disjoint trees.
\end{remark}

The concept of strong disjointness is the key ingredient behind the phenomenon of almost orthogonality responsible for the following  Bessel type inequalities.

\begin{proposition}[Nonmaximal Bessel's inequality, see \cite{Lac1}]
\label{Besselsineq1}
Let $j\in\{1,2\}$.
Consider a collection $\S'\subseteq\S$ of multitiles and let $\size_j(\cdot)$ denote the $j$-size  with respect to some function $F\in L^2(\R)$. Then $\S'$ can be written as a disjoint union 
$$\S'=\S_1'\cup\S_2',$$ 
where 
$$\size_j(\S_1')\le \frac{\size_j(\S')}{2}$$
while $\S_2'$ consists of a family $\F_{\S_2'}$ of pairwise disjoint trees  satisfying 
\begin{equation}
\label{e:intp}
\sum_{\T\in\F_{\S_2'}}|I_T|\lesssim \size_j(\S')^{-2}\|F\|_2^2,
\end{equation} 
with the implicit constant independent of $\S'$ and $F$. 
\end{proposition}

For $j=3$ we will need the following version of the above. We will call forest any collection $\F$ of strongly 3-disjoint trees and will denote by
$$N_{\F}(x)=\sum_{\T\in\F}1_{I_T}$$
the counting function of $\F$.

\begin{proposition}[Maximal Bessel's inequality]
\label{Besselsineq}
Let $\S'\subseteq\S$ be a collection  of multitiles and let $\U:=u_1<u_2<\ldots<u_J$ be an arbitrary sequence of integers. For each $s\in\S'$ let $j(s)$  de the unique number in $\{1,2,\ldots,J-1\}$ such that $2^{u_{j(s)}}\le |I_s|<2^{u_{j(s)+1}}$. Consider also an arbitrary sequence of functions $h_1,h_2,\ldots,h_{J-1}:\R\to\C$ satisfying 
$$\sum_{j=1}^{J-1}|h_j|^2\equiv 1,$$
and a function $h\in X_2(E)$, where $E$ is an arbitrary set of finite measure. 
Let $\size_3(\cdot)$ denote the 3-size   with respect to the functions $H_s:\R\to\C$ defined by 
$$H_s=hh_{j(s)}.$$ 
Then $\S'$ can be written as a disjoint union
$$\S'=\S_1'\cup\S_2',$$ 
where 
$$\size_3(\S_1')\le \frac{\size_3(\S')}{2}$$
while $\S_2'$ consists of a family $\F_{\S_2'}$ of pairwise disjoint trees  satisfying 
\begin{equation}
\label{e:intp33}
\sum_{\T\in\F_{\S_2'}}|I_T|\lesssim J^{1/8}\size_3(\S')^{-2}\left(\frac{1}{\size_3(\S')|E|^{1/2}}\right) ^{\epsilon}
\end{equation} 
for each $\epsilon>0$, with the implicit constant depending only on $\epsilon$.
\end{proposition}  

This proposition will follow from a chain of successive reductions, as in \cite{DTT}.

\begin{proposition}[Maximal Bessel's inequality, first reduction]
\label{Besselsineqred1}

Assume  $\S'\subseteq\S$ can be organized as a forest $\F'$ of trees $\T$ with tops $T$. 
Let $\U:=u_1<u_2<\ldots<u_J$ be an arbitrary sequence of integers. Consider also an arbitrary sequence of functions $h_1,h_2,\ldots,h_{J-1}:\R\to\C$ satisfying 
$$\sum_{j=1}^{J-1}|h_j|^2\equiv 1,$$
and a function $h\in X_2(E)$, where $E$ is an arbitrary set of finite measure. Define the functions $H_s$ as before. Assume also that
$$2^{m}\le \left(\frac1{|I_T|}\sum_{s\in\T}|\langle H_{s},\phi_{s,3}\rangle|^2\right)^{1/2}\le 2^{m+1}$$
for each $\T\in\F'$, and  
$$\left(\frac1{|I_{T'}|}\sum_{s\in\T\atop{I_s\subseteq I_{T'}}}|\langle H_{s},\phi_{s,3}\rangle|^2\right)^{1/2}\le 2^{m+1}$$
for each $T'\in\T\in\F'.$

Then 
$$\sum_{\T\in\F'}|I_T|\lesssim J^{1/8}2^{-2m}\left(\frac{1}{2^m|E|^{1/2}}\right)^{\epsilon}$$
for each $\epsilon>0$, with the implicit constant depending only on $\epsilon$.
\end{proposition}

Proposition \ref{Besselsineqred1} implies  Proposition \ref{Besselsineq} by a standard tree selection algorithm, see \cite{Lac1}.

The next reduction allows to replace the dependency on  $|E|^{1/2}$ with a dependency on the counting function multiplicity. By linearity we can also eliminate the size parameter $2^m.$

\begin{proposition}[Maximal Bessel's inequality, second reduction]
\label{Besselsineqred2}
Assume $\S'\subseteq\S$ can be organized as a forest $\F'$ of trees $\T$ with tops $T.$
Let $\U:=u_1<u_2<\ldots<u_J$ be an arbitrary sequence of integers. Consider also an arbitrary sequence of functions $h_1,h_2,\ldots,h_{J-1}:\R\to\C$ satisfying \begin{equation}
\label{lala44reqw3}
\sum_{j=1}^{J-1}|h_j|^2\equiv 1,
\end{equation}
and a function $h\in L^2(\R)$. Define the functions $H_s$ as before. Assume also that
\begin{equation}
\label{lala44reqw1}
1\le \left(\frac1{|I_T|}\sum_{s\in\T}|\langle H_{s},\phi_{s,3}\rangle|^2\right)^{1/2}\le 2
\end{equation}
for each $\T\in\F'$, and  
\begin{equation}
\label{lala44reqw2}
\left(\frac1{|I_{T'}|}\sum_{s\in\T\atop{I_s\subseteq I_{T'}}}|\langle H_{s},\phi_{s,3}\rangle|^2\right)^{1/2}\le 2
\end{equation}
for each $T'\in\T\in\F'.$ Let $I_0$ be an interval which contains the support of $N_{\F'}$.

Then 
\begin{equation}
\label{eq:drac}
\sum_{\T\in\F'}|I_T|\lesssim J^{1/8}\|N_{\F'}\|_{\infty}^{\epsilon}\int |h|^2\chi_{I_0}^{10}
\end{equation}
for each $\epsilon>0$, with the implicit constant depending only on $\epsilon$.
\end{proposition}

To see how Proposition \ref{Besselsineqred2} implies  Proposition \ref{Besselsineqred1} we first introduce the BMO norm of a forest $\F$ as
$$ \| \F\|_{\BMO} := \sup_I \frac{1}{|I|} \sum_{\T \in \F: I_T \subseteq I} |I_T|,$$ where the supremum is taken over all the 
dyadic intervals $I$. We then recall the following result from \cite{DTT}.

\begin{lemma}
\label{lem:goodlred}
Let $\F$ be a forest such that for some $\epsilon<1$ and some $A,B>0$
$$ \| N_{\F'}\|_1 \le A  \| N_{\F'}\|_{\infty}^{\epsilon} \hbox{ and } \| \F'\|_{\BMO} \le B  \|N_{\F'}\|_{\infty}^{\epsilon} $$
for all forests $\F' \subseteq \F$.  Then  we have
$$ \| N_{\F} \|_1 \lesssim_\epsilon A  B^{\frac{\epsilon}{1-\epsilon}}.$$
\end{lemma}

\begin{proof}[of Proposition \ref{Besselsineqred1} assuming Proposition \ref{Besselsineqred2}]  Let $\F' \subset \F$ be arbitrary.  From Proposition \ref{Besselsineqred2} with $h$ replaced by $h/2^m$, and $I_0$ chosen to be so large as to contain all the time intervals arising from $\F'$, we have
\begin{align*}
\| N_{\F'}\|_1 &= \sum_{\T \in \F'} |I_{T}| \\
&\lesssim
\sum_{s \in \bigcup_{\T \in \F'} \T}|\langle H_s/2^m,\phi_{s,3}\rangle|^2 \\
&\lesssim_\epsilon J^{1/8}\|N_{\F'}\|_{\infty}^{\epsilon} \int |h/2^m|^2 \\
&\lesssim J^{1/8} 2^{-2m} \|N_{\F'}\|_{\infty}^{\epsilon} 
\end{align*}
thanks to the $L^2$ normalization of $h \in X_2(E)$.  

On the other hand, if  $I_0$ is an arbitrary dyadic interval, then by
replacing $\F'$ with $\{ \T \in \F': I_T \subseteq I_0 \}$ in the above argument we see that
\begin{align*}
\frac{1}{|I_0|} \sum_{\T \in \F': I_T \subseteq I_0} |I_{\T}|
&\lesssim_\epsilon\frac{J^{1/8}}{|I_0|} \|N_{\F'}\|_{\infty}^{\epsilon} \int |f/2^m|^2 \chi_{I_0}^{10} \\
&\lesssim J^{1/8}\|N_{\F'}\|_{\infty}^{\epsilon} 2^{-2m} |E|^{-1} 
\end{align*}
thanks to the uniform bound of $|E|^{-1/2}$ on $f \in X_2(E)$.  Taking suprema over $I_0$ we conclude
that $\| \F' \|_{\BMO} \lesssim_\epsilon J^{1/8}\|N_{\F'}\|_{\infty}^{\epsilon} 2^{-2m} |E|^{-1}$.
Applying Lemma \ref{lem:goodlred} we conclude that
$$ \sum_{\T \in \F} |I_T| = \| N_{\F} \|_1 \lesssim_\epsilon J^{1/8}2^{-2m} (2^{-2m} |E|^{-1})^{\frac{\epsilon}{1-\epsilon}}.$$
Given the fact that $\epsilon>0$ is arbitrary, the proof of  Proposition \ref{Besselsineqred1} follows.
\end{proof}

We next focus on proving  Proposition \ref{Besselsineqred2}. 
We will borrow some of the terminology from \cite{DTT} in order to quote some results from there.
Let ${\mathcal D}_0, {\mathcal D}_1, {\mathcal D}_2$ be the dyadic grids
\begin{equation}\label{dyad}
\begin{split}
{\mathcal D}_0 &:= \{ [ 2^{-i} l, 2^{-i}(l+1)]: i, l \in \Z \}\;\;\hbox{(i.e. the standard dyadic grid)}   \\
{\mathcal D}_1 &:= \{ [ 2^{-i} (l + (-1)^i/3), 2^{-i}(l+1 + (-1)^i/3)]: i, l \in \Z \} \\
{\mathcal D}_2 &:= \{ [ 2^{-i} (l - (-1)^i/3), 2^{-i}(l+1 - (-1)^i/3)]: i, l \in \Z \}.
\end{split}
\end{equation}

One can easily verify that for every interval $J$ (not necessarily dyadic) there exists a $d \in \{0,1,2\}$ and a shifted 
dyadic interval $J' \in {\mathcal D}_d$ such that $J \subseteq J' \subseteq 3J$; we will say that $J$ is \emph{$d$-regular}.

Let $A \geq 1$, and let $d \in \{0,1,2\}$.  We shall say that a collection $\I \subset {\mathcal D}_0$ of time intervals 
is \emph{$(A,d)$-sparse} if we have the following properties:
\begin{itemize}
\item[(i)] If $I, I' \in \I$ are such that $|I| > |I'|$, then $|I| \geq 2^{100A} |I'|$.
\item[(ii)] If $I, I' \in \I$ are such that $|I| = |I'|$ and $I \neq I'$, then $\dist(I, I') \geq 100A |I'|$.
\item[(iii)] If $I \in \I$, then $A I$ is $d$-regular, thus there exists an interval $I_{A} \in {\mathcal D}_d$ such that
$AI \subseteq I_{A} \subseteq 3AI$.  We refer to $I_A$ as the \emph{$A$-enlargement} of $I$.
\end{itemize}

If $\I$ is an $(A,d)$-sparse set of time intervals and $P$ is a tile whose time interval $I_P$ lies in $\I$, we write $I_{P,A}$ for the
$A$-enlargement of $I_P$.  
We recall the following two results from \cite{DTT}.

\begin{theorem}
\label{thm:3b} Let $A, D> 1$ and
let $\F'$ be a forest with $\|N_\F'\|_{\infty} \leq D$.  Let also $\S' := \bigcup_{\T \in \F'} \T$, and suppose that the time intervals
$$ \{ I_s: s \in \S' \} \cup \{ I_T: \T \in \F' \}$$
are $(A,d)$-sparse.

Then there exists an exceptional set $\S_* \subset \S'$ of multitiles with
\begin{equation}\label{omega-bound-A}
|\bigcup_{s \in \S_*} I_s| \lesssim_\nu (A^{-\nu} + D^{-\nu}) \sum_{\T \in \F'} |I_T|
\end{equation}
such that we have the Bessel-type inequality
$$
\sum_{s \in \S' \backslash \S_*}|\langle f,\phi_{s,3}\rangle|^2 \lesssim_\nu 
((\log(2 + AD))^{10} + A^{10-\nu} D^{10}) \|f\|_{2}^2
$$
for each $\nu>1$ and each $f\in L^2(\R).$
\end{theorem}

\begin{lemma}[Sparsification] 
Let $\I$ be a collection of time intervals.  Then we can split $\I = \I_1 \cup \ldots \cup \I_L$ with $L = O(A^2)$ such that
each $\I_l$ for $1 \leq l \leq L$ is $(A,d)$-sparse for some $d=0,1,2$.
\end{lemma}

\begin{proof}[of Proposition \ref{Besselsineqred2}]
We start by noting that it suffices to prove  Proposition \ref{Besselsineqred2} without the localizing weight $\chi_{I_0}^{10}$. This is because $\chi_{I_0}^{-10}$ is a polynomial and hence $\chi_{I_0}^{-10}\psi_{s,3}$ satisfies the same properties ~\eqref{e.osc23} and ~\eqref{e.osc24} as $\psi_{s,3}$.
Let $\mu\ge 20$ be arbitrary.
We first apply the above lemma with $A=C_{\mu}(J\|N_{\F'}\|_{\infty})^{1/\mu}$ to the collection $\I:=\{I_s:s\in\S'\}$ to split $\I$ into $\I_1,\ldots,\I_{L}$, for some $L=O(A^2)$.
Each tree  $\T\in\F'$ will be disintegrated over the collections $\I_l$.
A further disintegration occurs by differentiating multitiles according to their spacial length, so in the end
$$\T=\bigcup_{1\le l\le L}\bigcup_{1\le l\le J-1}\S_{\T,l,j}.$$
Here $$\S_{\T,l,j}=\{s\in \T: I_s\in \I_l,\;2^{u_j}\le |I_s|<2^{u_{j+1}}\}$$
 
According to the Remark \ref{rem:md123}, each collection 
$$\F_{l,j}=\bigcup_{\T\in\F'}\bigcup_{\T'\in D(\S_{\T,l,j})}\T'$$
is a forest. It is easy to see that each $\F_{l,j}$ satisfies the requirements of Theorem \ref{thm:3b} with $A$ as above, and moreover
$$\|N_{\F_{l,j}}\|_{\infty}\le \|N_{\F'}\|_{\infty}$$
$$\|N_{\F_{l,j}}\|_{1}\le \|N_{\F'}\|_{1}.$$
Define
$$\S_{l,j}=\bigcup_{\T\in\F'}\S_{\T,l,j}.$$

By applying  Theorem \ref{thm:3b} to each forest $\F_{l,j}$ with $D=C_{\mu}J\|N\_{\F'}|_{\infty}$, $f=hh_j$ and $\nu=200\mu$ we get an exceptional set $\S_{l,j,*}$ such that  
\begin{equation}
\label{omega-bound-Aadded}
|\bigcup_{s \in \S_{l,j,*}} I_s| \lesssim_\mu C_{\mu}^{-200\mu}J^{-3}\|N_{\F'}\|_{\infty}^{-3} \sum_{\T \in \F_{l,j}} |I_T|
\end{equation}
and  we have the Bessel-type inequality
\begin{equation}
\label{sumnowneedfgrt}
\sum_{s \in \S_{l,j} \backslash \S_{l,j,*}}|\langle hh_j,\phi_{s,3}\rangle|^2 \lesssim_\mu 
((\log(2 + \|N_{\F'}\|_{\infty}))^{10}+(\log(2 + J))^{10}) \|hh_j\|_{2}^2.
\end{equation}
Define now 
$$\S_{*}:=\bigcup_{l,j}\S_{l,j,*}$$ and note that if $C_{\mu}$ is chosen sufficiently large then
\begin{equation}
\label{nufti6671}
|\bigcup_{s\in\S_{*}}I_s|\lesssim \frac{1}{10}\frac{\|N_{\F'}\|_{1}}{ \|N_{\F'}\|_{\infty}}.
\end{equation}
By summing up in \eqref{sumnowneedfgrt} and invoking \eqref{lala44reqw3} we get
\begin{align*}
\sum_{s \in \S' \backslash \S_{*}}|\langle H_s,\phi_{s,3}\rangle|^2 &\lesssim_\mu 
\|N_{\F'}\|_{\infty}^{\frac{2}{\mu}}J^{\frac{2}{\mu}}((\log(2 + \|N_{\F'}\|_{\infty}))^{10}+(\log(2 + J))^{10}) \|h\|_{2}^2\\&
\lesssim_\mu \|N_{\F'}\|_{\infty}^{\frac{3}{\mu}}J^{1/8}\|h\|_{2}^2.
\end{align*}

Write $\Omega := \bigcup_{s \in \S_*} I_s$.  Then
$$ \sum_{s \in\S': I_s \not \subseteq \Omega}|\langle H_s,\phi_{s,3}\rangle|^2\le \sum_{s \in \S' \backslash \S_*}|\langle H_s,\phi_{s,3}\rangle|^2.$$ 
To prove \eqref{eq:drac}, it thus suffices  to show that
\begin{equation}
\label{nufti6671hgg5}
 \sum_{s \in \S': I_s \subseteq \Omega}|\langle H_S,\phi_{s,3}\rangle|^2 \leq \frac{1}{2}
\sum_{s\in \S'}|\langle H_s,\phi_{s,3}\rangle|^2.
\end{equation}
From \eqref{lala44reqw1}, it thus suffices to show that
$$ \sum_{s\in\S': I_s \subseteq \Omega}|\langle H_s,\phi_{s,3}\rangle|^2 \leq \frac{1}{2} \| N_{\F'} \|_{1}.$$
For each tree $\T$ in $\F'$,
consider the multitile set $\{ s \in \T: I_s \subseteq \Omega \}$.  If $s$ is any tile in this set with $I_s$ maximal with respect to set inclusion, then $I_s \subseteq \Omega$ and from \eqref{lala44reqw2} we have
$$ \sum_{s' \in \T: I_{s'} \subseteq I_s}|\langle H_{s'},\phi_{s',3}\rangle|^2 \leq 4 |I_s|.$$
Summing this over all such $s$ upon noting that the $I_s$ are disjoint by dyadicity and maximality, we conclude that
$$ \sum_{s \in \T: I_s \subseteq \Omega}|\langle H_s,\phi_{s,3}\rangle|^2 \leq 4 |I_T \cap \Omega| = 4 \int_\Omega 1_{I_T}.$$
Summing this over all $\T \in \F'$ we obtain
$$ \sum_{s \in \S': I_s \subseteq \Omega}|\langle H_s,\phi_{s,3}\rangle|^2 \leq 4 \int_\Omega N_{\F'}
\leq 4 |\Omega| \| N_{\F'} \|_{\infty}$$
and the claim \eqref{nufti6671hgg5} follows now from \eqref{nufti6671}.
\end{proof}

\section{Single tree estimate}
\label{sec:single}

Consider a $3$-lacunary tree $\T$ and some coefficients $(c_s)_{s\in\T}$. The following representation holds for each $2^k\le |I_T|$, assuming $\Delta$ is chosen large enough 
\begin{equation}
\label{e:40}
\sum_{s\in \T:|I_s|\ge 2^k}c_{s}\psi_{s,3}(x)=\int\F(\sum_{s\in \T}c_{s}\psi_{s,3})(\xi)\zeta(2^k (\xi-c(\omega_{T,3})))e^{2\pi i\xi x}d\xi.
\end{equation}
Here $\zeta$ is some universal function equal to 1 on $[-100e,100e]$ and equal to 0 outside  $[-200e,200e]$. 

\begin{theorem}[Single tree estimate]
\label{thm:singletreeestimate}
Let $\T$ be a $3$-lacunary tree in $\S$ with top $T$, and let $$\\P_{\T}=\{I\:\:\text{dyadic}:I_s\subseteq I\subseteq I_{s'}\:\:\text{for some}\:\: s,s'\in \T\}$$ be the time convexification of $\T$.
Consider  a finite sequence of integers $u_1<u_2<\ldots<u_L$. For each $s\in\T$ let $l(s)$  be the unique number in $\{1,2,\ldots,L-1\}$ such that $2^{u_{l(s)}}\le |I_s|<2^{u_{l(s)+1}}$. Consider also an arbitrary sequence of functions $h_1,h_2,\ldots,h_{L-1}:\R\to\C$ satisfying 
$$\sum_{l=1}^{L-1}|h_l|^2\equiv 1,$$
and a function $h\in X(E)$, for some $E\subset \R$ of finite measure. For each $s\in\T$ define $H_s=hh_{l(s)}$.

Then
$$\left(\frac{1}{|I_T|}\sum_{s\in \T}|\langle H_s,\phi_{s,3}\rangle|^2\right)^{\frac12}\lesssim \sup_{I\in\\P_{\T}}\frac{1}{|I|}\int_{E} \chi_{I}^2,$$ with some universal implicit constant.
\end{theorem}
\begin{proof}

By frequency translation invariance we may assume that $0\in\omega_{T,3}$. If $h$ is supported outside $2I_T$ then from the decay of $\phi_{s,3}$ we get 
$$|\langle H_s,\phi_{s,3}\rangle|\lesssim\left(\frac{|I_s|}{|I_T|}\right)^{10}|I_T|^{-\frac{1}{2}}\int_{E} \chi_{I_s}^2.$$ This together with crude estimates based on the triangle inequality immediately prove the Theorem in this case.

 By invoking duality it hence suffices to prove for all $h$ supported on $2I_T$ and  all $(a_s)_{s\in\T}$ with $\|(a_s)\|_{l^2(\T)}\le 1$ that
$$\frac{1}{|I_T|^{\frac12}}\int h(\sum_{l=1}^{L-1}|\sum_{s\in\T\atop_{2^{u_l}\le |I_s|<2^{u_{l+1}}}}a_s\phi_{s,3}|^2)^{1/2}\lesssim \sup_{I\in\\P_{\T}}\frac{1}{|I|}\int_{E} \chi_{I}^2.$$ From ~\eqref{e:40} we know that 
\begin{align*}
|\sum_{s\in\T\atop_{2^{u_l}\le |I_s|<2^{u_{l+1}}}}a_s\phi_{s,3}(x)|&=\sup_{u_l\le k< u_{l+1}}|\sum_{s\in\T\atop{2^{k}\le |I_s|< 2^{u_{l+1}}}}a_s\psi_{s,3}(x)|\\&\lesssim M(\sum_{s\in\T\atop_{2^{u_l}\le |I_s|<2^{u_{l+1}}}}a_s\psi_{s,3}(x)).
\end{align*}
Define $$F_l(x)=\sum_{s\in\T\atop_{2^{u_l}\le |I_s|<2^{u_{l+1}}}}a_s\psi_{s,3}(x).$$
We will prove that
$$\frac{1}{|I_T|^{\frac12}}\int h(\sum_{l=1}^{L-1}M(F_l)^2)^{1/2}\lesssim \sup_{I\in\\P_{\T}}\frac{1}{|I|}\int_{E} \chi_{I}^2.$$ 

For a dyadic interval $J$ denote by $J_1,J_2,J_3$ the three dyadic intervals of the same length with $J$, sitting at the left of $J$, with $J_3$ being adjacent to $J$. Similarly  let $J_5,J_6,J_7$ be the three dyadic intervals of the same length with $J$, sitting at the right of $J$, with $J_5$ being adjacent to $J$. Also define $J_4=J$. Let $\mathcal J$  be the set of all dyadic intervals $J$ with the following properties:
\begin{enumerate}
\item[(a)]$ J\cap 2I_T\not=\emptyset$
\item[(b)]$ \nexists \;I\in\\P_{\T}: |I|<|J|\:\:and\:\: I\subset 3J$
\item[(c)]$ J_i\in \\P_{\T}\:\: for\:\: some\:\:1\le i\le 7.$ 
\end{enumerate}

We claim that $2I_T\subset\cup_{J\in\mathcal J}J$. Indeed, assume by contradiction that there exists some $x\in 2I_T\setminus\cup_{J\in\mathcal J}J$. Let $J^{(0)}\subset J^{(1)}\subset J^{(2)}\subset\ldots$ be the sequence of dyadic intervals of consecutive lengths containing $x$, with $|J^{(0)}|=\min_{I\in\\P_{\T}}|I|.$ Since $J^{(0)}\notin\mathcal J$ and since (a) and (b) are certainly satisfied for  $J^{(0)}$, it follows that $J_i^{(0)}\notin \\P_{\T}$ for each $1\le i\le 7$. Moreover, note that for each $1\le i\le 7$ there is no $I\in\\P_{\T}$ with $I\subset J_i^{(0)}$. We proceed now by induction. Assume that for some $j\ge 0$ we proved that for each $1\le i\le 7$ we have  $J_i^{(j)}\notin \\P_{\T}$  and also that there is no $I\in\\P_{\T}$ with $I\subset J_i^{(j)}$. Note that this implies the same for $j+1$. Indeed, since $3J^{(j+1)}\subset 7J^{(j)}$ and by induction hypothesis, it follows that (b) is satisfied for $J^{(j+1)}$. Hence $J_i^{(j+1)}\notin\\P_{\T}$ for each $1\le i\le 7$. We verify now  the second statement of the induction. Note that  if there was an $I\in\\P_{\T}$ with $I\subset J_i^{(j+1)}$ than the hypothesis of the induction and the fact that $3J^{(j+1)}\subset 7J^{(j)}$ would imply that $i\in\{1,2,6,7\}$. Hence $I\subset J_i^{(j+1)}\subset I_T$, and by convexity of $\\P_{\T}$ it would follow that $J_i^{(j+1)}\in\\P_{\T}$, impossible. This closes the induction. To see how the claim follows from here, observe that $I_T=J_i^{(j)}$ for some $i,j$, which certainly contradicts the fact that $I_T\in\\P_{\T}$. 

Next thing we prove is that on each interval $2J$ with $J\in \mathcal J$, the oscillation of $(\sum_{l=1}^{L-1}|F_j|^2)^{1/2}$ is well controlled. More exactly we will show that 
$$(\sum_{l=1}^{L-1}\sup_{x,y\in 2J}|F_l(x)-F_l(y)|^2)^{1/2}\lesssim \frac{1}{|J|^{\frac12}}.$$
Indeed, for each $x,y\in 2J$ we use  ~\eqref{e.osc23} and the fact that $0\in\omega_{T,3}$ to get 
\begin{align*}
(\sum_{l=1}^{L-1}\sup_{x,y\in 2J}|F_l(x)-F_l(y)|^2)^{1/2}&\lesssim |J|\sum_{s\in\T}\sup_{x\in 2J}|\frac{d}{dx}\psi_{s,3}(x)\|a_s|\\&\lesssim |J|\sum_{s\in\T}\frac{1}{|I_s|^{\frac32}}\chi_{I_s}^{10}(c(J))|a_s|,
\end{align*}
since by definition there exists no $s\in\T$ such that $I_s\subset 3J$.
Now 
\begin{align*}
|J|\sum_{s\in\T\atop{|I_s|\ge |J|}}\frac{1}{|I_s|^{\frac32}}\chi_{I_s}^{10}(c(J))|a_s|&\lesssim |J|(\sum_{s\in\T\atop{|I_s|\ge |J|}}|a_s|^2)^{\frac12}(\sum_{2^k\ge |J|}\sum_{i\gtrsim 1}\frac{1}{i^22^{3k}})^{\frac12}\\&\lesssim\frac1{|J|^{\frac12}}
\end{align*}
and also
\begin{align*}
|J|\sum_{s\in\T\atop{|I_s|<|J|}}\frac{1}{|I_s|^{\frac32}}\chi_{I_s}^{10}(c(J))|a_s|&\lesssim (\sum_{s\in\T\atop{|I_s|< |J|}}|a_s|^2)^{\frac12}(\sum_{2^k<|J|}\sum_{i\gtrsim \frac{|J|}{2^k}}\frac{1}{i^22^{3k}})^{\frac12}\\&\lesssim\frac1{|J|^{\frac12}}, 
\end{align*}
due to the fact that  there exists no $I_s\subset 3J$ with $|I_s|<|J|$.

Define now the measure space $X=\cup_{J\in\mathcal J}J$ and its sigma-algebra $\Upsilon$ generated by the maximal intervals $ J\in\mathcal J$. Recall that $2I_T\subset\cup_{J\in\mathcal J}J=X\subset 10I_T$. We will see that for each $x\in J$ 
\begin{equation}
\label{eq:50}
(\sum_{l=1}^{L-1}M(F_l)^2(x))^{1/2}\lesssim \frac1{|J|}\int_{J} (\sum_{l=1}^{L-1}M(F_l)^2(z))^{1/2}dz+\frac{1}{|J|^{\frac12}}.
\end{equation}
Indeed, if $r_l>\frac12|J|$ then  
\begin{align*}
(\sum_{l=l}^{L-1}(\frac{1}{2r_j}\int_{x-r_l}^{x+r_l}|F_l|(z)dz)^2)^{1/2}&\lesssim \inf_{y\in J}(\sum_{l=1}^{L-1}M(F_l)^2(y))^{1/2}\\&\lesssim \frac1{|J|}\int_{J}(\sum_{l=1}^{L-1}M(F_l)^2(z))^{1/2}dz
.
\end{align*}
 On the other hand, if $r_l\le\frac12|J|$ then
\begin{align*}
(\sum_{l=1}^{L-1}(\frac{1}{2r_l}\int_{x-r_l}^{x+r_l}|F_l|(z)dz)^2)^{1/2}&\lesssim \sup_{y\in 2J}(\sum_{l=1}^{L-1}|F_l|^2(y))^{1/2}\\&\lesssim \inf_{y\in J}(\sum_{l=1}^{L-1}|F_l|^2(y))^{1/2}+\frac{1}{|J|^{\frac12}}\\&\lesssim \frac1{|J|}\int_{J}(\sum_{l=1}^{L-1}M(F_l)^2(z))^{1/2}dz+\frac{1}{|J|^{\frac12}}.
\end{align*}
From ~\eqref{eq:50} we can write
\begin{align*}
\frac{1}{|I_T|^{\frac {1}{2}}}\int h(\sum_{l=1}^{L-1}M(F_l)^2)^{1/2}&\lesssim \frac{1}{|I_T|^{\frac {1}{2}}}\int_{X} h\E((\sum_{l=1}^{L-1}M(F_l)^2)^{1/2}|\Upsilon)+\sup_{J\in \mathcal J}\frac {1}{|J|}\int_{J}|h|
\\&=\frac{1}{|I_T|^{\frac12}}\int_{X} \E(h|\Upsilon)\E((\sum_{l=1}^{L-1}M(F_l)^2)^{1/2}|\Upsilon)+\sup_{J\in \mathcal J}\frac1{|J|}\int_{J}|h|\\&\le \frac{1}{|I_T|^{\frac {1}{2}}}\|\E(h|\Upsilon)\|_{L^{\infty}}\int_{X} \E((\sum_{l=1}^{L-1}M(F_l)^2)^{1/2}|\Upsilon)+
\sup_{J\in  \mathcal J}\frac{1}{|J|}\int_{E} \chi_{J}^{2}\\&\lesssim \left(1+\left[\int_{X} \E((\sum_{l=1}^{L-1}M(F_l)^2)^{1/2}|\Upsilon)^2\right]^{\frac12}\right)\sup_{J\in  \mathcal J}\frac{1}{|J|}\int_{E} \chi_{J}^{2}\\&\lesssim \sup_{J\in  \mathcal J}\frac{1}{|J|}\int_{E} \chi_{J}^{2},
\end{align*}
where $\E(\cdot|\Upsilon)$ denotes the conditional expectation relative to $\Upsilon.$ Finally, note that since for each $J\in\mathcal J$, $J_i\in\\P_{\T}$ for some $i$, we have that 
$$\sup_{J\in  \mathcal J}\frac{1}{|J|}\int_{E} \chi_{J}^2\lesssim \sup_{I\in \\P_{\T}}\frac{1}{|I|}\int_{E} \chi_{I}^2,$$ which certainly ends the proof of our theorem.
\end{proof}

\section{The proof of Theorem ~\ref{t.oscprincmodel}}
\label{proof}
For each $f,g\in L^2(\R)$ and each subset of tiles $\S'\subset\S$ define
$$M_{\S'}(f,g)(x)=(\sum_{j=1}^{J-1}|\sum_{s\in\S'\atop{2^{u_j}\le|I_s|<2^{u_{j+1}}}}|I_s|^{-1/2}\langle f,\phi_{s,1}\rangle\langle g,\phi_{s,2}\rangle\phi_{s,3}|^2)^{1/2}.$$
In order to prove inequality ~\eqref{e.form2} it suffices to show that
$$|\{x:M_{\S}(f,g)(x)>\lambda\}|\lesssim \frac{J^{1/4}}{\lambda},$$
uniformly for each $f,g$ with $\|f\|_2=\|g\|_2=1$ and each $\lambda>0.$
We claim that it suffices to prove the above for $\lambda\sim 1.$ Indeed, for arbitrary $\lambda$ choose $k\in\Z$ such that $2^{k}\le \lambda<2^{k+1}$. Also, for each multitile 
$$s=[m2^i,(m+1)2^i]\times\prod_{j=1}^{3}[\frac{l_j}{5}2^{-i},(\frac{l_j}{5}+1)2^{-i}]\in\S$$ define $$s(k):=[m2^{i+k},(m+1)2^{i+k}]\times\prod_{j=1}^{3}[\frac{l_j}{5}2^{-i-k},(\frac{l_j}{5}+1)2^{-i-k}]\in\S.$$
Define also the collection of multitiles
$$\S(k)=\{s(k):s\in\S\},$$
and the functions
$$\psi_{s(k),j}(x)=\frac{1}{2^{k/2}}\psi_{s,j}(\frac{x}{2^k}).$$
Note that $\S(k)$ and $\psi_{s(k),j}$ satisfy all the requirements of Theorem ~\ref{t.oscprincmodel} (with a different choice of grids), with the same implicit constants in ~\eqref{e.osc23}. The claim now follows by noting that

$$|I_s|^{-1/2}\langle f,\phi_{s,1}\rangle\langle g,\phi_{s,2}\rangle \phi_{s,3}(x)=|I_{s(k)}|^{-1/2}\langle \operatorname{Dil}_{2^k}^2 f,\phi_{s(k),1}\rangle\langle \operatorname{Dil}_{2^k}^2 g,\phi_{s(k),2}\rangle2^k\phi_{s(k),3}(2^kx).$$

Define now  $\Gamma=\max\{[-\log_2(\size_1(\S))],[-\log_2(\size_2(\S))]\}$, where the 1-size is understood here with respect to the function $f$ while the 2-size is understood with respect to $g$. 
By iterating Proposition ~\ref{Besselsineq1} simultaneously for $F=f$ and $F=g$, it follows that $\S$ can be written as a disjoint union $\S=\bigcup_{n\ge \Gamma}\S_n,$ with 
\begin{equation}
\label{gftyrrcv1}
\size_j(\bigcup_{m\ge n}\S_m)\le 2^{-n}
\end{equation} 
for $j=1,2$, and each $\S_n$ consists of a family $\F_{\S_n}$ of pairwise disjoint trees  satisfying 
\begin{equation}
\label{gftyrrcv2}
\sum_{\T\in\F_{\S_n}}|I_T|\lesssim 2^{2n}.
\end{equation}
The contributions coming from the collections  $\S'=\bigcup_{\Gamma\le n\le 0}\S_n$ and $\S''=\bigcup_{n>0}\S_n$ will be evaluated quite differently.

In the case of  $\S'$, crude estimates will suffice. Let $\T\in\F_{\S_n}$. By using ~\eqref{gftyrrcv1}, the decay in ~\eqref{e.osc23} and the triangle inequality we immediately get the estimate
\begin{equation}
\label{gftyrrcv3}
M_{\T}(f,g)(x)\lesssim 2^{-2n}\chi_{I_T}^{4}(x)
\end{equation}
for each $x\notin 2I_T$. For each $\Gamma\le n\le 0$ and each $\T\in\F_{\S_n}$ set $E_n=2^{-n}I_T$. From ~\eqref{gftyrrcv2} we get that the exceptional set $E=\bigcup_{\Gamma\le n\le 0}E_n$ has measure $\lesssim 1.$ Also, ~\eqref{gftyrrcv3} implies that
\begin{align*}
\|M_{\S'}(f,g)\|_{L_1(E^c)}&\le \sum_{\Gamma\le n\le 0}\|M_{\S_n}(f,g)\|_{L_1(E_n^c)}\\&\lesssim \sum_{\Gamma\le n\le 0}2^{-2n}2^{3n}\\&\lesssim 1.
\end{align*}
We conclude that 
\begin{equation}
\label{gftyrrcv4}
|\{x:M_{\S'}(f,g)(x)\gtrsim 1\}|\lesssim 1.
\end{equation}

We will focus next on the estimates for the collection $\S''$. This time we will rely on the fact that $\size_1(\S'')\le 1$ and $\size_2(\S'')\le 1.$ 
As before, for each $s\in\S$ let $j(s)$  denote the unique number in $\{1,2,\ldots,J-1\}$ such that $2^{u_{j(s)}}\le |I_s|<2^{u_{j(s)+1}}$.
The 3-size will now intervene in a crucial way. Define
$$V:=\{x:M_{\S''}(f,g)(x)\gtrsim 1\}.$$ If $|V|\le 1$ there is nothing to prove, so we will assume $|V|>1$.
Let $h_1,h_2,\ldots,h_{J-1}:\R\to\C$
 be functions satisfying
$$\sum_{j=1}^{J-1}|h_j|^2\equiv 1$$ such that
$$M_{\S''}(f,g)(x)=
\sum_{s\in\S''}|I_s|^{-1/2}\langle f,\phi_{s,1}\rangle\langle g,\phi_{s,2}\rangle h_{j(s)}\phi_{s,3}(x).
$$

From Theorem ~\ref{thm:singletreeestimate} we know that the 3-size of the collection $\S''$ with respect to the functions $H_s:=|V|^{-\frac12}1_{V}h_{j(s)}$ is $\lesssim 1.$ There is actually no restriction in assuming it is $\le 1$. By applying iteratively the propositions ~\ref{Besselsineq1} and ~\ref{Besselsineq} to the collection $\S''$  it follows that $\S''$ can be written as a disjoint union $\S''=\bigcup_{n\ge 0}\S_n'',$ with 
\begin{equation}
\label{gftyrrcv1222}
\size_j(\S_n'')\le 2^{-n}
\end{equation} 
for $j\in\{1,2,3\}$, and each $\S_n''$ consists of a family $\F_{\S_n''}$ of pairwise disjoint trees  satisfying 
\begin{equation}
\label{gftyrrcv32222}
\sum_{\T\in\F_{\S_n''}}|I_T|\lesssim J^{1/8}2^{\frac52n}.
\end{equation}
Finally, by Lemma ~\ref{l:328}
\begin{align*}
|V|^{1/2}&\lesssim  \langle M_{\S''}(f,g),|V|^{-\frac12}1_{V}\rangle\\&=\sum_{s\in\S''}|I_s|^{-1/2}\langle f,\phi_{s,1}\rangle\langle g,\phi_{s,2}\rangle\langle H_s, \phi_{s,3}\rangle\\&\le \sum_{n\ge 0}\sum_{s\in\S_n''}|I_s|^{-1/2}\langle f,\phi_{s,1}\rangle\langle g,\phi_{s,2}\rangle\langle H_s, \phi_{s,3}\rangle\\&\lesssim \sum_{n\ge 0}2^{-3n}\sum_{\T\in\F_{\S_n''}}|I_T|\\&\lesssim \sum_{n\ge 0}2^{-3n}2^{\frac52n}J^{1/8}\\&\lesssim J^{1/8} .
\end{align*}

We conclude that 
\begin{equation}
\label{gftyrrcv5}
|\{x:M_{\S''}(f,g)(x)\gtrsim 1\}|\lesssim J^{1/4}.
\end{equation}
The estimates in ~\eqref{gftyrrcv4} and ~\eqref{gftyrrcv5} end the proof of Theorem ~\ref{t.oscprincmodel}.

\section{The proof of Theorem ~\ref{t.square}}
\label{sec:square}
By applying discretization techniques like in Section \ref{sec:dis}, Theorem~\ref{t.square} will follow from the following. 
\begin{theorem}
\label{t.oscsecmodel}
Let $\G$, $\G_1$, $\G_2$ and $\G_3$ be four grids with $\G$ satisfying
\begin{equation*}
\label{e.thmgrid1'}
\G\subset\B
\end{equation*}
\begin{equation*}
\label{e.thmgrid2'}
I,I'\in\G\Rightarrow\max\{|I||I'|^{-1},|I'||I|^{-1}\}\ge 2^{\Delta}. 
\end{equation*}
Let $e$ be a  number with $10^2<|e|\le 10^5$ and define
$$\D=\{\prod_{j=1}^{3}[\frac{l_j}{5}2^i,(\frac{l_j}{5}+1)2^i]\in\prod_{j=1}^{3}\G_j:\; l_2=l_1+e,\;l_3=l_1+l_2\}.$$
Define also the  set of multitiles 
$$\S=\{s=I_{s}\times Q_{s}: I_s\in\G,\:\: Q_s\in\D\;\text{with sidelength}\;\frac{1}{|I_s|}\}.$$
Assume that  each multitile $s=I_s\times\prod_{j=1}^{3}\omega_{s,3}\in\S$ is associated with three functions $(\psi_{s,j})_{j=1}^{3}$ satisfying
\begin{equation}
\label{e.osc23'}
|\frac{d^n}{d x^n}\operatorname{Mod}_{-c(\omega_{s,j})}\psi_{s,j}(x)|\lesssim_{n,M} |I_s|^{-1/2-n}\chi_{I_s}^M,\;\;n,M\ge 0
\end{equation}
\begin{equation}
\label{e.osc24'}
\operatorname{supp}\widehat{\psi}_{s,j}\subset \omega_{s,j},
\end{equation}
for each $j=1,2,3.$ 

Then  for each $f,g\in L^2(\R)$ we have the estimate
\begin{equation}
\label{e.form12222}
\|\sup_{k\in\Z}|\sum_{s\in\S\atop{|I_s|=2^k}}|I_s|^{-1/2}\langle f,\psi_{s,1}\rangle\langle g,\psi_{s,2}\rangle\psi_{s,3}|^2)^{1/2}\|_{1,\infty}\lesssim \|f\|_2\|g\|_2,
\end{equation}
with the implicit constant depending only on the implicit constants in ~\eqref{e.osc23'}.
\end{theorem}

To prove this theorem amounts to proving Theorem ~\ref{t.oscprincmodel} in the case $u_1,u_2,\ldots,u_{J}$ are consecutive integers, with a bound independent of $J$. By analyzing the whole argument for  Theorem ~\ref{t.oscprincmodel}, it follows that the dependency on $J$ in there comes from a single source, that is the Maximal Bessel inequality in Proposition \ref{Besselsineq}. This dependency is eliminated by proving the following version of Proposition \ref{Besselsineqred1}.

\begin{proposition}
\label{Besselsineqred1'}

Assume  $\S'\subseteq\S$ can be organized as a forest $\F'$ of trees $\T$ with tops $T$.  Consider  an arbitrary sequence of functions $h_k:\R\to\C$, $k\in\Z$ satisfying 
$$\sum_{k\in\Z}|h_k|^2\equiv 1,$$
and a function $h\in L^{2}(\R)$. Define the functions $H_s$ by $H_s=hh_{\log(|I_s|)}$. Assume also that
$$2^{m}\le \left(\frac1{|I_T|}\sum_{s\in\T}|\langle H_{s},\psi_{s,3}\rangle|^2\right)^{1/2}$$
for each $\T\in\F'$.

Then 
$$\sum_{\T\in\F'}|I_T|\lesssim 2^{-2m}\|h\|_2^2.$$
\end{proposition}
\begin{proof}
It suffices to prove that
$$\sum_{s\in\S'\atop{|I_s|=2^k}}|\langle f,\psi_{s,3}\rangle|^2\lesssim \|f\|_2^2,$$
uniformly in all $k\in\Z$ and all $f\in L^2(\R)$. This in turn will follow by duality from the following estimate
$$\|\sum_{s\in\S'\atop{|I_s|=2^k}}a_s\psi_{s,3}\|_2^2\lesssim \sum_{s\in\S'\atop{|I_s|=2^k}}|a_s|^2,$$
which holds for all sequence $(a_s)$.
\end{proof}
Indeed,
\begin{align*}
\|\sum_{s\in\S'\atop{|I_s|=2^k}}a_s\psi_{s,3}\|_2^2&=\sum_{\omega\in\G_3:|\omega|=2^{-k}}\|\sum_{s\in\S'\atop{\omega_s=\omega}}a_s\psi_{s,3}\|_2^2\\&=\sum_{\omega\in\G_3:|\omega|=2^{-k}}\sum_{s,s'\in\S'\atop{\omega_s=\omega_{s'}=\omega}}a_s\bar{a}_{s'}\langle\psi_{s,3},\psi_{s',3}\rangle\\&\lesssim \sum_{\omega\in\G_3:|\omega|=2^{-k}}\sum_{s,s'\in\S'\atop{\omega_s=\omega_{s'}=\omega}}|a_sa_{s'}|(1+\frac{\operatorname{dist}(I_s,I_{s'})}{2^k})^{-10}\\&\lesssim \sum_{\omega\in\G_3:|\omega|=2^{-k}}\sum_{s\in\S'\atop{\omega_s=\omega}}|a_s|^2\\&=\sum_{s\in\S'}|a_s|^2.
\end{align*}

\end{document}